\documentclass[11pt,twoside,leqno]{aomamlt2e} 
  \pageno{911}
\received{February 20, 2004}

\newcommand\odinger{Schr{\hskip2pt\rm \"{\hskip-7.5pt\it o}}dinger}
\newcommand\hfq{\hfill\qed}

\let\newpf\proof \let\proof\relax

\def\be{\begin{equation}}

\def\op{\overline \partial}

\def\ba{{\begin{align}}}
\def\ea{{\end{align}}}

\def\Ad{{\mathrm {Ad}}}
\def\pa{\partial}

\def\g{{\gamma}}

\def\SO{{\mathrm {SO}}}
\def\GL{{\mathrm {GL}}}

\def\Sl{\mathrm {sl}}

\def\H{{\mathbb H}}

\def\rot{\operatorname{rot}}

\def\0{{\mathbf 0}}

\def\cal{\mathcal}

\newcommand{\Sr}{{\mathbb S}}

\newcommand{\Diff}{\mathrm{Diff}}

\newtheorem{thm}{Theorem}[section]
 
\newtheorem{cor}[thm]{Corollary}

\newtheorem{lemma}[thm]{Lemma}

\newtheorem{prop}[thm]{Proposition}

\theorembodyfont{\upshape}
\newtheorem{rem}{{\it Remark}}[section]

\numberwithin{equation}{section}

\def\theequation {\thesection.\arabic{equation}}

\def \bn {\hfill \\ \smallskip\noindent}

\theoremstyle{definition}

\def\proof{\bn {\bf Proof.} }

\def\note#1
{\marginpar
{\tiny $\leftarrow$
\par
\hfuzz=20pt \hbadness=9000 \hyphenpenalty=-100 \exhyphenpenalty=-100
\pretolerance=-1 \tolerance=9999 \doublehyphendemerits=-100000
\finalhyphendemerits=-100000 \baselineskip=6pt
#1}\hfuzz=1pt}

\newcommand{\id}{\operatorname{id}}

\newcommand{\MM}{{\cal M}}

\newcommand{\C}{{\mathbb C}}

\newcommand{\Q}{{\mathbb Q}}
\newcommand{\R}{{\mathbb R}}
\newcommand{\T}{{\mathbb T}}

\newcommand{\Z}{{\mathbb Z}}

\def\B0{{\bold{0}}}

\catcode`\@=12

\def\Empty{}
\newcommand\oplabel[1]{
  \def\OpArg{#1} \ifx \OpArg\Empty {} \else
  	\label{#1}
  \fi}

\newcommand{\comm}[1]{}
\newcommand{\comment}[1]{}
 \def\twoauthors#1#2{\vskip15pt
\hbox to\hsize{\hss%
\small By {\scshape#1}\ifacks\global\acksfalse*\fi\   and  {\scshape#2}\hss}
\ifshort\else
\newcommand{\theauthors}{{\scriptsize\scshape\uppercase{#1 and #2}}}\fi
\vskip7pt
\vskip\baselineskip
\global\authortrue\everypar={\global\authorfalse\everypar={}}}

\begin{document}
\currannalsline{164}{2006} 

 \title{Reducibility or nonuniform hyperbolicity\\ for quasiperiodic
Schr\"odinger cocycles}

 \acknowledgements{A. A. is a Clay Research Fellow.}
%\author{}
\twoauthors{Artur Avila}{Rapha\"el Krikorian}

 \institution{Coll\`ege de France, Paris, France\\
\current{CNRS UMR 7599, Universit\'e Pierre et Marie Curie,\\ Paris, France}\\
\email{artur@ccr.jussieu.fr}\\
\vglue-9pt
Universit\'e Pierre et Marie Curie, Paris, France\\
\email{krikoria@ccr.jussieu.fr}}

 \shorttitle{Reducible or nonuniformly hyperbolic Schr\"odinger cocycles}

\centerline{\bf Abstract}
\vglue6pt

We show that for almost every frequency $\alpha \in \R \setminus \Q$,
for every $C^\omega$ potential $v:\R/\Z \to \R$, and for almost every energy
$E$ the corresponding quasiperiodic Schr\"odinger cocycle
is either reducible or nonuniformly hyperbolic.
%(similar results are valid in the smooth category).
This result gives   very good control on
the absolutely continuous part of the spectrum of the corresponding
quasiperiodic Schr\"odinger operator, and allows us to complete the
proof of the Aubry-Andr\'e conjecture on the measure of the spectrum of the
Almost Mathieu Operator.

\section{Introduction}

A {\it one-dimensional quasiperiodic $C^r$-cocycle in ${\rm SL}(2,\R)$}
(briefly, a $C^r$-co\-cycle) is a
pair $(\alpha,A) \in \R \times C^r(\R/\Z,{\rm SL}(2,\R))$,
viewed as a {\it linear skew-product}\/:
\begin{align} \label {skpro}
(\alpha,A):\R/\Z \times \R^2 &\to \R/\Z \times \R^2\\
\nonumber
(x,w) &\mapsto (x+\alpha,A(x) \cdot w).
\end{align}

For $n \in \Z$, we let $A_n \in C^r(\R/\Z,{\rm SL}(2,\R))$ be defined by the
rule $(\alpha,A)^n=(n\alpha,A_n)$ (we will keep
the dependence of $A_n$ on $\alpha$ implicit).  Thus $A_0(x)=\id$,
\begin{equation}
A_n(x)=\prod_{j=n-1}^0 A(x+j\alpha)=A(x+(n-1)\alpha) \cdots A(x), \quad
\text {for } n \geq 1,
\end{equation}
and $A_{-n}(x)=A_n(x-n\alpha)^{-1}$.
The Lyapunov exponent of $(\alpha,A)$ is defined as
\begin{equation}
L(\alpha,A)=\lim_{n \to \infty}
\frac {1} {n} \int_{\R/\Z} \ln \|A_n(x)\| dx \geq 0.
\end{equation}
%where $A_n(x)=\prod_{j=n-1}^0 A(x+j\alpha)=A(x+(n-1)\alpha) \cdots A(x)$
%(we will keep the dependence on $\alpha$ implicit).

Also,  $(\alpha,A)$ is {\it uniformly hyperbolic} if there exists
a continuous splitting $E_s(x) \oplus E_u(x)=\R^2$, and $C>0$,
$0<\lambda<1$ such that for every $n\geq 1$ we have
\begin{align}
\|A_n(x) \cdot w\| \leq C \lambda^n \|w\|, \quad w \in E_s(x),\\
\nonumber
\|A_{-n}(x) \cdot w\| \leq C \lambda^n \|w\|, \quad w \in E_u(x).
\end{align}
Such splitting is automatically unique and thus invariant; that is,
$A(x)E_s(x)=E_s(x+\alpha)$ and $A(x)E_u(x)=E_u(x+\alpha)$.
The set of uniformly hyperbolic cocycles is open in the $C^0$-topology
(one allows perturbations both in $\alpha$ and in $A$).

Uniformly hyperbolic cocycles have a positive Lyapunov exponent.
If $(\alpha,A)$ has positive Lyapunov exponent but is not uniformly
hyperbolic then it will be called {\it nonuniformly hyperbolic}.

We say that a $C^r$-cocycle $(\alpha,A)$ is {\it $C^r$-reducible}
if there exists $$B \in  C^r(\R/2\Z,{\rm SL}(2,\R)) 
\hbox{ and }A_* \in {\rm SL}(2,\R)$$ such that
\begin{equation} \label {reducib}
B(x+\alpha) A(x) B(x)^{-1}=A_*, \quad x \in \R.
\end{equation}
Also,  $(\alpha,A)$ is $C^r$-reducible modulo $\Z$ if one can take $B
\in C^r(\R/\Z,{\rm SL}(2,\R))$.\footnote{Obviously, reducibility modulo $\Z$ is a
stronger notion than plain reducibility, but in some situations one can show
that both definitions are equivalent (see Remark~1.5).
The advantage of defining reducibility ``modulo $2\Z$'' is to include some
special situations (notably certain uniformly hyperbolic cocycles).}

Now, $\alpha \in \R \setminus \Q$ satisfies a {\it Diophantine
condition} ${\rm DC}(\kappa,\tau)$, $\kappa>0$, $\tau>0$ if
\begin{equation} \label {dc}
\left |q \alpha-p\right |>\kappa |q|^{-\tau}, \quad (p,q) \in
\Z^2,\quad q \neq 0.
\end{equation}
Let ${\rm DC}=\cup_{\kappa>0,\tau>0} {\rm DC}(\kappa,\tau)$.  It is well known that
$\cup_{\kappa>0} {\rm DC}(\kappa,\tau)$ has full Lebesgue measure if $\tau>1$.

Now,  $\alpha \in \R \setminus \Q$ satisfies a {\it recurrent
Diophantine condition} ${\rm RDC}(\kappa,\tau)$ if there are
infinitely many $n>0$ such that $G^n(\{\alpha\}) \in {\rm DC}(\kappa,\tau)$, where
$\{\alpha\}$ is the fractional part of $\alpha$ and $G:(0,1) \to [0,1)$
is the Gauss map $G(x)=\{x^{-1}\}$.  We let
${\rm RDC}=\cup_{\kappa>0,\tau>0} {\rm RDC}(\kappa,\tau)$.  Notice that
${\rm RDC}(\kappa,\tau)$ has full Lebesgue measure as long as ${\rm DC}(\kappa,\tau)$
has positive Lebesgue measure (since the Gauss map is
ergodic with respect to the probability measure $\frac {dx} {(1+x) \ln 2}$).
It is possible to show that $\R \setminus {\rm RDC}$ has Hausdorff dimension $1/2$.

%\marginpar {Check}

Given $v \in C^r(\R/\Z,\R)$, let us consider the Schr\"odinger cocycle
\begin{equation}
S_{v,E}(x)=
\begin{pmatrix}
E-v(x) & -1\\
1 & 0
\end{pmatrix} \in C^r(\R/\Z,{\rm SL}(2,\R))
\end{equation}
($v$ is called the potential and $E$ is called the energy).

There is  fairly good comprehension of the dynamics of Schr\"odinger
cocycles in the case of either small or large potentials:

\begin{prop}[Sorets-Spencer \cite {SS}] \label {sorets-spencer}
Let $v \in C^\omega(\R/\Z,\R)$ be a nonconstant potential{\rm ,} and let $\alpha
\in \R$.  There exists $\lambda_0=\lambda_0(v)>0$ such that if
$|\lambda|>\lambda_0$ then for every $E \in \R$ there is
$L(\alpha,S_{\lambda v,E})>0$.
\end{prop}

\setcounter{thm}{2}
 {\scshape Proposition 1.2} (Eliasson [E1]\footnote{This result was originally
stated for the continuous time case, but the proof also works for the
discrete time case.}).
{\it Let $v \in C^\omega(\R/\Z,\R)${\rm ,} and let $\alpha \in {\rm DC}$.  There exists
$\lambda_0=\lambda_0(v,\alpha)$ such that if $|\lambda|<\lambda_0$ then for
almost every $E \in \R$ the cocycle $(\alpha,S_{\lambda v,E})$ is
$C^\omega$-reducible.}
\vskip4pt

  {\it Remark} 1.1.
Sorets-Spencer's result is {\it nonperturbative}: the
``largeness'' condition $\lambda_0$ does not depend on $\alpha$.
On the other hand, the proof of Eliasson's result is
{\it perturbative}: the ``smallness'' condition
$\lambda_0$ depends in principle on $\alpha$ (in the full measure set ${\rm DC}
\subset \R$).  We will come back to this issue (cf.\ Theorem~\ref{nonperturbative}).

\vskip4pt {\it Remark} 1.2.
In general, one cannot replace ``almost every'' by ``every''
in Eliasson's result above.  Indeed, in \cite {eliasson} it is also shown
that the set of energies for which
$(\alpha,S_{\lambda v,E})$ is not (even $C^0$)
reducible is nonempty for a generic (in an appropriate topology)
choice of $(\lambda,v)$ satisfying $|\lambda|<\lambda_0(v)$.
Those ``exceptional'' energies do have zero Lyapunov exponent.

\vskip4pt {\it Remark} 1.3. % \label {unifhyperbol}
Let $\alpha \in {\rm DC}$ and $A \in C^r(\R/\Z,{\rm SL}(2,\R))$, $r=\infty,\omega$.
In this case, $(\alpha,A)$ is uniformly hyperbolic if and only if it is
$C^r$-reducible and has a positive Lyapunov exponent, see \cite[\S 2]{E2}.
Thus, there are lots of ``simple cocycles'' for which one has
positive Lyapunov exponent, resp.\ reducibility, and indeed both at the same
time: this is the case in particular for $|E|$ large in the Schr\"odinger
case.  Those examples are also stable (here we fix $\alpha \in {\rm DC}$ and
stability is with respect to perturbations of $A$).
 
However, cocycles with a positive Lyapunov exponent, resp.\ reducible,
but which are {\it not} uniformly hyperbolic do happen for a
{\it positive measure set of energies}
for many choices of the potential, and in particular in the situations
described by the results of Sorets-Spencer (this follows from
\cite[Th.~12.14]{B}), resp.\ Eliasson.
\vskip4pt
\setcounter{rem}{5}

Our main result for Schr\"odinger cocycles aims to close the gap and
describe the situation (for almost every energy) without
largeness/smallness assumption on the potential:

\vskip4pt {\scshape Theorem A.} {\it Let $\alpha \in {\rm RDC}$ and let $v:\R/\Z \to \R$ be a $C^\omega$
potential. Then{\rm ,} for Lebesgue almost every $E${\rm ,} the cocycle
$(\alpha,S_{v,E})$ is either
nonuniformly hyperbolic or $C^\omega$-reducible.}
\vskip4pt 

For $\theta \in \R$, let
\begin{equation}
R_\theta=
\begin{pmatrix}
\cos 2\pi\theta & -\sin 2\pi\theta\\
\sin 2\pi\theta & \cos 2\pi\theta
\end{pmatrix}.
\end{equation}
Given a $C^r$-cocycle $(\alpha,A)$, we associate a canonical
one-parameter family of $C^r$-cocycles $\theta \mapsto
(\alpha,R_\theta A)$.
Our proof of Theorem~A goes through for the more general context of
cocycles homotopic to the identity, with the role of the energy parameter
replaced by the $\theta$ parameter.

\vskip6pt {\scshape Theorem A$'$.}
{\it Let $\alpha \in {\rm RDC}${\rm ,} and let
$A:\R/\Z \to {\rm SL}(2,\R)$ be $C^\omega$
and homotopic to the identity\/\footnote{For the case of cocycles
nonhomotopic to the identity, see \cite {AK}.}.
Then for Lebesgue almost every $\theta \in \R/\Z${\rm ,} the cocycle
$(\alpha,R_\theta A)$ is either nonuniformly hyperbolic or
$C^\omega$-reducible.}

\vglue6pt {\it Remark} 1.4. 
Theorems~A and~A$'$ also hold in the smooth setting.  The only modification in
the proof is in the use of a KAM theoretical result of Eliasson
(see Theorem~\ref {eliassoncocycle}), which must be replaced by a smooth
version.  They also generalize to the case of continuous time (differential
equations): in this case the adaptation is straightforward.  See \cite {AK3}
for a discussion of those generalizations.
%and their applications to the study of the
%quasiperiodic Schr\"odinger operator.

\pagegoal=50pc
\vglue6pt {\it Remark} 1.5.  
One can distinguish two distinct behaviors among the
reducible cocycles $(\alpha,A)$ given by Theorems A and A$'$.  The first
is uniformly hyperbolic behavior; see Remark~1.3.
The second is {\it totally elliptic} behavior, corresponding (projectively)
to an irrational rotation of $\T^2 \equiv \R/\Z \times \mathbb{P}^1$.
More precisely, we call a cocycle totally elliptic if it is $C^r$-reducible
and the constant matrix $A_*$ in (\ref {reducib}) can be chosen to be a
rotation $R_\rho$, where $(1,\alpha,\rho)$ are linearly independent
over $\Q$.  In this case it is easy to see that
the cocycle $(\alpha,A)$ is automatically $C^r$-reducible modulo $\Z$
(possibly replacing $\rho$ by $\rho+\frac {\alpha} {2}$).  (To see that
almost every reducible cocycle is either uniformly hyperbolic or totally
elliptic, it is enough to use Theorems \ref {dioph} and \ref {dioph1} which
are due to Johnson-Moser and Deift-Simon.)
\vglue6pt

Theorems A and A$'$ give a nice global picture for the theory of
quasiperiodic cocycles, extending known results for cocycles taking values
on certain {\it compact} groups (see \cite {K1} for the case of
$\mathrm {{\rm SU}}(2)$).  They fit with the Palis conjecture for
general dynamical systems \cite {Pa}, and have a strong analogy with
the work of Lyubich in the quadratic family \cite {Ly}, generalized in
\cite {ALM}.

More importantly, reducible and nonuniformly
hyperbolic systems can be efficiently described through a wide variety of
methods, especially in the analytic case.  With respect to reducible
systems, the dynamics of the cocycle itself is of course very simple, and
the use of KAM theoretical methods (\cite {DS}, \cite {eliasson}) allowed
also a good comprehension of their perturbations.  With respect to
nonuniformly hyperbolic systems, there has
been recently lots of success in the application of
subtle properties of subharmonic functions (\cite {BG}, \cite {GS}, \cite
{BJ}) to obtain large deviation estimates with important consequences
(such as regularity properties of the Lyapunov exponent).

\Subsec{Application to \odinger\ operators}
We now discuss the application of the previous results to
the quasiperiodic Schr\"odinger operator
\begin{equation}
H_{v,\alpha,x} u(n)=u(n+1)+u(n-1)+v(x+\alpha n) u(n), \quad u \in l^2(\Z),
\end{equation}
where $\alpha \in \R \setminus \Q$, $x \in \R$ and
$v:\R/\Z \to \R$ is $C^\omega$.
The properties of $H_{v,\alpha,x}$ are closely
connected to the properties of the family of
cocycles $(\alpha,S_{v,E})$, $E \in \R$.  Notice for instance that
if $(u_n)_{n \in \Z}$ is a solution of $H_{v,\alpha,x}u=Eu$ then
\begin{equation}
\begin{pmatrix}E-v(x+n\alpha)&-1\\1&0\end{pmatrix} \cdot
\begin{pmatrix}u(n)\\u(n-1)\end{pmatrix}=
\begin{pmatrix}u(n+1)\\u(n)\end{pmatrix}.
\end{equation}

Let $\Sigma$ be the spectrum of $H_{v,\alpha,x}$.  It is
well known (see \cite {JM})
%for a proof in the continuous time case
that
\begin{equation}
\Sigma=\{E \in \R,\, (\alpha,S_{v,E}) \text { is {\it not} uniformly
hyperbolic}\},
\end{equation}
so that $\Sigma=\Sigma(v,\alpha)$ does not depend on $x$.

Let $\Sigma_{sc}=\Sigma_{sc}(\alpha,v,x)$ (respectively, $\Sigma_{ac}$,
$\Sigma_{pp}$) be (the support of) the singular continuous
(respectively, absolutely continuous, pure point) part of the spectrum
of $H_{v,\alpha,x}$.

It has been shown by Last-Simon (\cite {LS}, Theorem 1.5) that
$\Sigma_{ac}$ does not depend on $x$ for $\alpha \in \R \setminus \Q$
(there are no hypotheses on the smoothness of $v$ beyond continuity).
It is known that
$\Sigma_{sc}$ and $\Sigma_{pp}$ do depend on $x$ in general.

We will also introduce some decompositions
of $\Sigma$ that only depend on the cocycle, and hence are independent of $x$.

We split $\Sigma=\Sigma_0 \cup \Sigma_+$ in the
parts corresponding to zero Lyapunov exponent and positive Lyapunov
exponent for the cocycle $(\alpha,S_{v,E})$.
By \cite {BJ}, $\Sigma_0$ is closed.
%By Proposition~\ref {bj continuity}, $\Sigma_0$ is closed.
% in the real analytic case.

Let $\Sigma_r$ be the set of $E \in \Sigma$ such that
$(\alpha,S_{v,E})$ is $C^\omega$-reducible.
%Let $\Sigma_{ar}$ be the set of $E \in \Sigma$ such that
%$(\alpha,S_{v,E})$ is almost $C^infty$-reducible.
%Using Proposition~\ref {eliasson
%open}, we see that $\Sigma \setminus \Sigma_{ar}$ is closed.
It is easy to see that $\Sigma_r \subset \Sigma_0$.

Notice that by the Ishii-Pastur Theorem (see \cite {I} and \cite {P}),
we have $\Sigma_{ac}\! \subset\! \Sigma_0$.
%$\Sigma_{ac} \subset \overline
%\Sigma_0' \subset \overline \Sigma_0$, where $\Sigma_0'$ is the set of
%density points of $\Sigma_0$ (more precisely, $\Sigma_0$ is an ``essential
%support'' for the absolutely continuous spectrum).

By Theorem A, $\Sigma_0 \setminus \Sigma_r$ has zero
Lebesgue measure if $\alpha \in {\rm RDC}$ and $v \in C^\omega$.
One way to interpret $|\Sigma_0 \setminus \Sigma_r|=0$ (using the
Ishii-Pastur Theorem) is that generalized
eigenfunctions in the essential support of the absolutely continuous
spectrum are (very regular) Bloch waves.  This already
gives (in the particular cases under consideration) strong
versions of some conjectures in the literature (see for instance the
discussion after Theorem 7.1 in \cite {DeiftSimon}).  (Analogous statements
hold in the continuous time case.)

Another immediate application of Theorem A is a nonperturbative version of
Eliasson's result stated in Proposition~1.2.  It is based on the
following nonperturbative result:

\begin{prop}[Bourgain-Jitomirskaya] \label {acbj}
Let $\alpha \in {\rm DC}${\rm ,} $v \in C^\omega$.  There exists
$\lambda_0=\lambda_0(v)>0$ \/{\rm (}\/only depending on the bounds of $v${\rm ,} but not on
$\alpha${\rm )} such that if $|\lambda|<\lambda_0${\rm ,} then
the spectrum of $H_{\lambda v,\alpha,x}$ is purely absolutely continuous
for almost every $x$.
\end{prop}

\begin{thm} \label {nonperturbative}
Let $\alpha \in {\rm RDC}${\rm ,} $v \in C^\omega$.  There exists
$\lambda_0>0$ \/{\rm (}\/which may be taken the same as in the previous
proposition\/{\rm )}\/ such that if $|\lambda|<\lambda_0${\rm ,} then
$(\alpha,S_{\lambda v,E})$ is reducible for almost every $E$.
\end{thm}

\Proof 
By the previous proposition, $\Sigma_{ac}=\Sigma$, so that $\Sigma_+=\emptyset$.
\Endproof 

There are several other interesting results which can be concluded easily
from Theorem~A and current results and techniques:

\smallbreak\noindent \hskip5pt   (1)
 Zero Lebesgue measure of $\Sigma_{sc}$ for almost every frequency,
\smallbreak\noindent \hskip5pt \hangindent=23pt\hangafter=1
(2) Persistence of absolutely continuous spectrum
under perturbations of the potential,
\smallbreak\noindent \hskip5pt \hangindent=23pt\hangafter=1
(3) Continuity of the Lebesgue measure of $\Sigma$ under perturbations of
the potential.

 \smallbreak\noindent 
Although the key ideas behind those results are quite
transparent (given the appropriate background), a proper
treatment would take us too far from the proof of Theorem~A, which is the
main goal of this paper.  We will thus concentrate on a particular case
which provides one of the most striking applications of Theorem~A.  For the
applications mentioned above (and others), see \cite {AK3}.

\vglue-5pt
\Subsubsec{Almost Mathieu}
Certainly the most studied family of potentials
in the literature is $v(\theta)=\lambda \cos 2 \pi \theta$, $\lambda>0$.
In this case, $H_{v,\alpha,x}$ is called the Almost Mathieu Operator.
%Before discussing applications for general potentials, it is
%worth to exemplify some consequences in this context.

The Aubry-Andr\'e conjecture on the measure of the spectrum of the Almost
Mathieu Operator states that the measure of the spectrum of $H_{\lambda
\cos 2\pi \theta,\alpha,x}$ is $|4-2\lambda|$ for every
$\alpha \in \R \setminus \Q$, $x \in \R$ (see \cite {AA}).\footnote
{The ``critical case'' $\lambda=2$ can be traced even further
back to Hofstadter \cite {H}.}
There is a long story of
developments around this problem, which led to several partial results
(\cite {HS}, \cite {AMS}, \cite {L}, \cite {JK}).
In particular, it has already been proved
for every $\lambda \neq 2$ (see \cite {JK}), and for every $\alpha$ not of
constant type\footnote {A number $\alpha \in \R$ is said to be of {\it
constant type} if the coefficients of its continued fraction expansion are
bounded.  It follows that $\alpha$ is of constant type if and only if
$\alpha \in \cup_{\kappa>0} {\rm DC}(\kappa,1)$ if and only if $\alpha \in
\cup_{\kappa>0} {\rm RDC}(\kappa,1)$.} \cite {L}.
However, for $\alpha$, say, the golden mean, and $\lambda=2$, where one
should prove zero Lebesgue measure of the spectrum,
previous to this work, it was  still unknown even whether the spectrum has empty interior.

Using Theorem~A, we can deal with the
last cases (which are also Problem~5 of \cite {Simon}).
%See also Simon, Schr\"odinger Operators in the Twentieth First Century,
%Problem 5.

\begin{thm} \label {almost}
\hskip-5pt The spectrum of $H_{\lambda \cos 2\pi \theta,\alpha,x}$ has Lebesgue measure
$|4\!-\!2\lambda|$ for every $\alpha \in \R \setminus \Q$.

\end{thm}

\Proof 
As stated above, it is enough to consider $\lambda=2$ and $\alpha$ of
constant type, in particular $\alpha \in {\rm RDC}$.
%\footnote {Let us note however
%that Theorem~A can be used to derive new proofs of other cases, such as
%$\lambda \neq 2$ and $\alpha$ of constant type.}
Let $\Sigma$ be the spectrum of $H_{2 \cos 2\pi\theta,\alpha,x}$.
%Consider the
%decomposition $\Sigma=\Sigma_r \cup \Sigma_+ \cup (\Sigma_0 \setminus
%\Sigma_r)$.
By Corollary~2 of \cite {BJ}, $\Sigma_+=\emptyset$.
By Theorem~A, for
almost every $E \in \Sigma_0$, $(\alpha,S_{2 \cos 2\pi\theta,E})$ is
$C^\omega$-reducible.  Thus, it is enough to show
that $(\alpha,S_{2 \cos 2\pi\theta,E})$ is not $C^\omega$-reducible for
every $E \in \Sigma$.

Assume this is not the case, that is, $(\alpha,S_{2 \cos 2\pi\theta,E})$ is
reducible for some $E \in \Sigma$. 
To reach a contradiction, we will approximate the potential
$2 \cos 2\pi\theta$ by $\lambda \cos 2\pi \theta$ with
$\lambda>2$ close to $2$.
Then, by Theorem A of \cite {eliasson}, if $(\lambda,E')$ is sufficiently
close to $(2,E)$, either $(\alpha,S_{\lambda \cos 2\pi\theta,E'})$ is
uniformly hyperbolic or $L(\alpha,S_{\lambda \cos 2\pi\theta,E'})=0$.
In particular
(since the spectrum depends continuously on the potential), there exists $E'
\in \R$ such that $L(\alpha,S_{\lambda \cos 2\pi\theta,E'})=0$.
But it is well known, see \cite {H}, that the Lyapunov exponent of
$S_{\lambda \cos 2\pi\theta,E'}$ is bounded from below by $\max
\{\ln \frac {\lambda}{2},0\}>0$ and the result follows.
%\footnote{Alternatively, we could have used \cite {eliasson} to conclude
%that there exists some absolutely continuous spectrum for $H(\lambda \cos
%\theta,\alpha,x)$, which contradicts positivity of the Lyapunov exponent
%(Ishii-Pastur).  We could also argue by keeping the potential constant and
%varying $\alpha$.  Indeed, any $\alpha$ of constant type can be
%approximated by $\alpha' \in {\rm DC}(\kappa,\tau)$ of nonconstant type
%(with $\kappa$, $\tau$ fixed).  For such
%$\alpha'$, absence of absolutely
%continuous spectrum follows from the zero measure of the spectrum proved in
%\cite {L}, and this gives a contradiction with \cite {eliasson} as
%before.}.
\hfq

\begin{rem}
Barry Simon has pointed out to us an alternative argument based on duality
that shows that if $\alpha \in \R \setminus \Q$ and if $E \in
\Sigma=\Sigma(2 \cos 2\pi\theta,\alpha)$ then the cocycle
$(\alpha,S_{2 \cos 2\pi\theta,E})$ is not $C^\omega$-reducible.  Indeed, if
$(\alpha,S_{v,E})$ is $C^\omega$-reducible and $E \in \Sigma$,
then (by duality) there exists $x \in \R$ such that $E$ is an eigenvalue
for $H_{2 \cos 2\pi\theta,\alpha,x}$, and the corresponding eigenvector
decays exponentially, hence $L(\alpha,S_{v,E})>0$ which gives a
contradiction.  (This argument actually can be used to show that
$(\alpha,S_{v,E})$ is not $C^1$-reducible.)
\end{rem}

By \cite {GJLS}, we get:

\begin{cor}
The spectrum of $H_{2\cos 2\pi\theta,\alpha,x}$ is purely
singular continuous for every $\alpha \in \R \setminus \Q${\rm ,} and for almost
every $x \in \R/\Z$.
\end{cor}

Theorem A also gives a fairly precise dynamical picture for $\lambda<2$
(completing the spectral picture obtained by Jitomirskaya in \cite {J}):

\begin{thm} \hskip-5pt
Let $\lambda\!<\!2${\rm ,} $\alpha\! \in\! {\rm RDC}$.  For almost every $E\! \in\! \R${\rm ,} 
$(\alpha,S_{\lambda \cos 2\pi\theta,E})$ is reducible.
\end{thm}

\Proof 
By Corollary 2 of \cite {BJ}, the Lyapunov exponent is zero on the spectrum.
%(that $|\Sigma_+|=0$ also follows from the work of Jitomirskaya \cite {J}
%and is enough for our purposes).
The result is now a consequence of Theorem A.
\hfq

\Subsec{Outline of the proof of Theorem~{\rm A}}
The proof has some distinct steps, and is based on a renormalization scheme.
This point of view, which has already been used  in the study of
reducibility properties of quasiperiodic cocycles with values in ${\rm SU}(2)$
and ${\rm SL}(2,\R)$, has proved to be very useful in the nonperturbative case
(see \cite{K1}, \cite{K2}).  However, the scheme we present in this paper
is somehow simpler and fits better (at least in the ${\rm SL}(2,\R)$ case)
with the general renormalization philosophy (see  \cite {S} for a very nice
description of this point of view on renormalization):

%The proof has some distinct steps, and is based on a renormalization scheme
%(see \cite {S} for a very nice description of this point of view on
%renormalization):

\begin{enumerate}

\item The starting point is the theory of Kotani\footnote{This step holds in
much greater generality, namely for cocycles over ergodic transformations.}.
For almost every energy $E$, if the Lyapunov exponent of $(\alpha,S_{v,E})$
is zero, then the cocycle is $L^2$-conjugate to a cocycle in
$\SO(2,\R)$.  Moreover, the {\it fibered rotation number} of the cocycle
%(which is closely related to the integrated density of states)
is Diophantine with respect to $\alpha$.  (The set $\Delta$ of those
energies will be precisely the set of energies for which we will be able to
conclude reducibility.)

\item We now consider a smooth
cocycle $(\alpha,A)$ which is $L^2$-conjugate to rotations.  An
explicit estimate allows us to control the derivatives of iterates of the
cocycle restricted to certain small intervals.

\item After introducing the notion of renormalization of cocycles,
we interpret item (2) as ``{\it a priori\/} bounds'' (or precompactness)
for a sequence of renormalizations $(\alpha_{n_k},A^{(n_k)})$.

\item The recurrent Diophantine condition for $\alpha$ allows us
to take $\alpha_{n_k}$ uniformly
Diophantine, so that the limits of renormalization are
cocycles $(\hat \alpha,\hat A)$ where $\hat \alpha$ satisfies a
Diophantine condition.  Those limits are essentially (that is, modulo a
constant conjugacy) cocycles in $\SO(2,\R)$,
and are trivial to analyze: they are always reducible.

\item Since $\lim (\alpha_{n_k},A^{(n_k)})$ is reducible, Eliasson's theorem
\cite {eliasson}
allows us to conclude that some renormalization
$(\alpha_{n_k},A^{(n_k)})$ must be reducible, provided the
fibered rotation number of
$(\alpha_{n_k},A^{(n_k)})$ is Diophantine with respect to $\alpha_{n_k}$.

\item This last condition is actually equivalent to the fibered rotation
number of $(\alpha,A)$ being Diophantine with respect to $\alpha$.
It is easy to see that reducibility is invariant under
renormalization and  so $(\alpha,A)$ is itself reducible.

%\item The Integrated Density of States $N(E)$
%can be identified with the fibered rotation number of $(\alpha,S_{v,E})$,
%so reducibility follows for the set $\Delta$
%of energies given by step (1).

\end{enumerate}

We conclude that for almost every $E \in \R$ such that
$L(\alpha,S_{v,E})=0$, the cocycle $(\alpha,S_{v,E})$ is reducible, which
is equivalent to Theorem A by Remark~1.3.

The above strategy uses $\alpha \in {\rm RDC}$ in order to take
good limits of renormalization.  It would be interesting to try to obtain
results under the weaker condition $\alpha \in {\rm DC}$ by working directly with
deep renormalizations (without considering limits).

\begin{rem}

Renormalization methods have been previously applied to the study of
quasiperiodic Schr\"odinger operators, see for instance \cite {BF}, \cite
{FK} and \cite {HS}.  While the
notions used by Helffer-Sj\"ostrand are quite different from ours, the
``monodromization techniques'' of Buslaev-Fedotov-Klopp
correspond to essentially the
same notion of renormalization used here.  An important conceptual
difference is in the use of renormalization: we are interested in the
dynamics of the renormalization operator itself, in a spirit close to
works in
one-dimensional dynamics (see for instance \cite {Ly}, \cite {Y},
\cite {S}).
\end{rem}

\section{Parameter exclusion}
\vglue-12pt
\Subsec{$L^2$-estimates}
We say that $(\alpha,A)$ is $L^2$-conjugated to a cocycle of rotations
if there exists a measurable
$B:\R/\Z \to {\rm SL}(2,\R)$ such that $\|B\| \in L^2$ and
\begin{equation}
B(x+\alpha) A(x) B(x)^{-1} \in {\rm SO}(2,\R).
\end{equation}

\begin{thm} \label {l2}
Let $v:\R/\Z \to \R$ be continuous.
Then for almost every $E${\rm ,} either $L(\alpha,S_{v,E})>0$
or $S_{v,E}$ is $L^2$-conjugated to a cocycle of rotations.
\end{thm}

\Proof 
Looking at the projectivized action of $(\alpha,S_{v,E})$ on the
upper half-plane $\H$, one sees that
the existence of an $L^2$ conjugacy to rotations is equivalent to
the existence of a measurable invariant section\footnote {That is
$S_{v,E}(x) \cdot m(x,E)=m(x+\alpha,E)$.} $m(\cdot,E):\R/\Z \to \H$
satisfying $\int_{\R/\Z} \frac {1} {\Im m(x,E)} dx<\infty$.  This holds
for almost every $E$ such that $L(\alpha,S_{v,E})=0$ by Kotani Theory, as
described in \cite {Simon1}\footnote {This reference was pointed out to us 
by Hakan Eliasson.} (the measurable invariant section $m$ we want is given
by $\frac {-1} {m_-}$ in the notation of \cite {Simon1}).
\Endproof\vskip4pt 

It turns out that this result generalizes to the
setting of Theorem A$'$:

\begin{thm} \label {l3}

Let $A:\R/\Z \to {\rm SL}(2,\R)$ be continuous.  Then for almost every $\theta
\in \R${\rm ,} either $L(\alpha,R_\theta A)>0$ or
$(\alpha,R_\theta A)$ is $L^2$-conjugated to a cocycle of rotations.
\end{thm}

The proof of this generalization is essentially the same as in the
Schr\"odinger case.  We point the reader to \cite {AK} for a discussion of
this and further generalizations.

\begin{rem}
Both theorems above are valid in a much more general setting, namely for
cocycles over transformations preserving a probability measure.  The
requirement on the cocycle is the least  to speak of Lyapunov exponents
(and Oseledets theory), namely integrability of the logarithm of the norm.
\end{rem}

\vglue-12pt
\Subsec{Fibered rotation number}
Besides the Lyapunov exponent, there is one important invariant
associated to continuous cocycles which are
homotopic to the identity.  This invariant, called the {\it fibered rotation
number} will be denoted by $\rho(\alpha,A) \in \R/\Z$, and
was introduced in \cite {H}, \cite {JM} (we recall its definition in
Appendix A).
The fibered rotation number is a continuous function of
$(\alpha,A)$, where $(\alpha,A)$ varies in the space of continuous cocycles
which are homotopic to the identity.  Another important elementary
fact is that both $E \mapsto -\rho(\alpha,S_{v,E})$ and $\theta \mapsto
\rho(\alpha,R_\theta A)$ have nondecreasing lifts $\R \to \R$, and in
particular, those functions have nonnegative derivatives almost everywhere.
The following result was proved in \cite {JM}, in the continuous time case,
and in \cite {DeiftSimon}, in the discrete time case used here
(and where an optimal estimate is given).

\begin{thm} \label {dioph}
Let $v \in C^0(\R/\Z,\R)$.  Then for almost every $E$ such that
$L(\alpha,S_{v,E})=0${\rm ,} 
\begin{equation}
\frac {d} {dE} \rho(\alpha,S_{v,E})<0.
\end{equation}
\end{thm}

This result (and proof) also generalize to the setting of Theorem A$'$ (see
\cite {AK} for further generalizations):

\begin{thm} \label {dioph1}
Let $A \in C^0(\R/\Z,{\rm SL}(2,\R))$ be continuous and homotopic to the
identity.  Then for almost every $E$ such that
$L(\alpha,R_\theta A)=0${\rm ,}  
\begin{equation} 
\frac {d} {d\theta}
\rho(\alpha,R_\theta A)>0.
\end{equation} 
\end{thm}

\begin{rem}
In the Schr\"odinger case, it is possible to show that the fibered rotation
number is a surjective function (of $E$) onto $[0,1/2]$.  In \cite {AS} it
is also shown that $N(E)=1-2\rho(\alpha,S_{v,E})$ can be interpreted as
the integrated density of states.
\end{rem}

The arithmetic properties of the fibered rotation number are also important
for the analysis of cocycles $(\alpha,A)$. 
Fix $\alpha \in \R$.
Let us say that $\beta \in \R/\Z$ is {\it Diophantine with
respect to $\alpha$} if there exists $\kappa>0$, $\tau>0$ such that
\begin{equation} \label {2.4}
\|2 \beta-k\alpha\|_{\R/\Z} \geq \kappa (1+|k|)^{-\tau}, \quad k \in \Z,
\end{equation}
where $\| \cdot \|_{\R/\Z}$ denotes the distance to the nearest integer.
If $\tau>1$ then the Lebesgue measure of the set of $\beta \in \R/\Z$ 
which satisfy (\ref {2.4}) is at least $1-2 \frac {\tau+1} {\tau-1} 
\kappa$.  In particular, Lebesgue almost every $\beta$ is Diophantine with 
respect to $\alpha$.  By Theorems 2.3 and 2.4 we conclude:

\begin{cor} \label {dioph2}

Let $\alpha \in {\rm DC}${\rm ,} $v \in C^0(\R/\Z,\R)$.  Then for almost every $E
\in \R$ such that $L(\alpha,S_{v,E})=0${\rm ,} 
$\rho(\alpha,S_{v,E})$ is
Diophantine with respect to $\alpha$.

\end{cor}

\begin{cor} \label {dioph3}

Let $\alpha \in {\rm DC}${\rm ,} $A \in C^0(\R/\Z,{\rm SL}(2,\R))$.  Then for almost every
$\theta \in \R$ such that $L(\alpha,R_\theta A)=0${\rm ,} 
$\rho(\alpha,R_\theta A)$ is Diophantine with respect to $\alpha$.
\end{cor}

\vglue-4pt
The fibered rotation number and its arithmetic properties play a role in the
following result of Eliasson \cite {eliasson}:

\begin{thm} \label {eliassoncocycle}
Let $(\alpha,A) \in \R \times C^\omega(\R/\Z,{\rm SL}(2,\R))$.  Assume that\/{\rm :}\/
 \smallbreak \noindent \hskip6pt
{\rm (1)} $\alpha \in {\rm DC}(\kappa,\tau)$ for some $\kappa>0${\rm ,} $\tau>0${\rm ,}
\smallbreak \noindent \hskip6pt
{\rm (2)} $\rho(\alpha,A)$ is Diophantine with respect to \pagebreak $\alpha${\rm ,}

\smallbreak \noindent \hskip6pt
{\rm (3)} $A$ admits a holomorphic extension to some strip $\R/\Z \times
(-\epsilon,\epsilon)${\rm ,}

\smallbreak \noindent \hskip6pt
{\rm (4)} $A$ is sufficiently close to a constant $\hat A \in {\rm SL}(2,\R)${\rm :}
\begin{equation}
\sup_{z \in \R/\Z \times (-\epsilon,\epsilon)} \|A(z)-\hat A\|<
\delta=\delta(\kappa,\tau,\epsilon,\hat A).
\end{equation}
Then $(\alpha,A)$ is reducible.
\end{thm}

This theorem was originally proved in the case of differential equations,
but the adaptation to our setting is immediate.  For further
generalizations, see \cite {AK3}.

\section{Estimates for derivatives}

%For the rest of this paper, we will go back to the setting of the
%introduction, that is, we will consider cocycles $(\alpha,A)$ where $\alpha
%\in \R \setminus \Q$ and $A:\R/\Z \to {\rm SL}(2,\R)$.
In this section, we will assume that $(\alpha,A)$ is $L^2$-conjugated
to a cocycle of rotations.  There exist measurable
$B:\R/\Z \to {\rm SL}(2,\R)$ and $R:\R/\Z\to \SO(2,\R)$ such that
\begin{equation}
\forall x\in\R/\Z, \quad
A(x)=B(x+\alpha)R(x)B(x)^{-1}\quad {\rm and}\quad
\int_{\R/\Z}\phi(x)dx<\infty
\end{equation}
where we set $\phi(x)=\|B(x)\|^2=\|B(x)^{-1}\|^2$ (here and in
what follows, $\R^2$ is supplied with the Euclidean norm and the space of
real $2 \times 2$ matrices $\mathrm{M}(2,\R)$
is supplied with the operator norm).
%\max \{\|B(x)\|^2,\|B(x)^{-1}\|^2\}=\|B(x)\|^2$

We introduce the {\it maximal function} $S(\cdot)$ of $\phi$: 
\begin{equation}
S(x)=\sup_{n \geq 1}\frac{1}{n}\sum_{k=0}^{n-1}\phi(x+k\alpha).
\end{equation}
Since the dynamics of $x\mapsto x+\alpha$ is ergodic on
$\R/\Z$ endowed with Lebesgue measure, the
Maximal Ergodic Theorem gives us the weak-type inequality
\begin{equation}
\forall M>0,\quad {\rm Leb}(\{x\in\R/\Z,S(x)>M\}) \leq
\frac{1}{M}\int_{\R/\Z}\phi(x)dx,
\end{equation}
and for a.e $x_0\in\R/\Z$ the quantity $S(x_0)$ is finite.
%Notice that if $S(x_0)<\infty$, then $\phi(x_0+k\alpha),S(x_0+k\alpha)
%\leq (k+1) S(x_0)$ are finite for all $k \geq 0$.

If $X \in \GL(2,\R)$, we let $\Ad(X)$ be the linear operator in
$\mathrm{M}(2,\R)$ which is given by
$\Ad(X) \cdot Y=X \cdot Y \cdot X^{-1}$.  Notice that the operator norm
of $\Ad(X)$ satisfies the bound $\|\Ad(X)\| \leq \|X\| \cdot \|X^{-1}\|$.

\begin{lemma}

Assume that $A$ is Lipschitz \/{\rm (}\/with constant ${\rm Lip}(A)${\rm ).}
%and $L^2$-conjugate
%to a cocycle of rotations by $B$.
Then for every $x_0,x \in \R/\Z$ such that $S(x_0)<\infty${\rm ,} 
\begin{equation} \label {lip}
\|A_n(x_0)^{-1}(A_n(x)-A_n(x_0))\|\leq
e^{n|x-x_0|\|A\|_{C^0}{\rm Lip}(A)\phi(x_0)S(x_0)}-1,
\end{equation}
and in particular
\begin{equation} \label {3.5}
\|A_n(x)\| \leq e^{n |x-x_0|\|A\|_{C^0}{\rm Lip}(A)S(x_0)\phi(x_0)}
\biggl (\phi(x_0)\phi(x_0+n\alpha) \biggr )^{1/2}.
\end{equation}
\end{lemma}

\Proof 
We compute $I_n(x_0,x):=A_n(x_0)^{-1}(A_n(x)-A_n(x_0))$:
\begin{align}
&\hskip-6pt I_n(x_0,x)\\
&=
A_n(x_0)^{-1} \biggl (\prod_{k=n-1}^0 \bigl
(A(x_0+k\alpha)+(A(x+k\alpha)-A(x_0+k\alpha)) \bigr )-A_n(x_0) \biggr )\nonumber\\[6pt]
\nonumber
&=\sum_{r=1}^n\sum_{0\leq i_r<\ldots<i_1\leq n-1}
\prod_{j=1}^r (\Ad(A_{i_j}(x_0)^{-1}) \cdot H_{i_j}(x_0,x))
\end{align}
where we have set 
\begin{equation}
H_i(x_0,x)=A(x_0+i\alpha)^{-1} \cdot (A(x+i\alpha)-A(x_0+i\alpha)),
\end{equation}
so that
\begin{equation}
\|H_i(x_0,x)\| \leq \|A\|_{C^0} {\rm Lip}(A) |x-x_0|.
\end{equation}
The assumptions we made give
\begin{equation}
\|A_i(x_0)\|=\|A_i(x_0)^{-1}\|\leq \|B(x_0+i\alpha)^{-1}\|
\cdot \|B(x_0)\|;
%\quad \|A_{i-1}(x_0)\|\leq \|B(x_0+i\alpha)\| \cdot \|B(x_0)^{-1}\|;
\end{equation}
that is,
\begin{equation}
\|\Ad(A_i(x_0)^{-1})\|\leq (\|B(x_0+i\alpha)^{-1}\| \cdot
\|B(x_0)\|)^2=\phi(x_0)\phi(x_0+i\alpha).
\end{equation}

Thus 
\begin{align}
\|I_n(x_0,x)\|&
\leq \sum_{r=1}^n\sum_{0\leq i_r<\ldots<i_1\leq n-1}
\prod_{j=1}^r \biggl (\|A\|_{C^0} {\rm Lip}(A) |x-x_0|
\phi(x_0) \phi(x_0+i_j \alpha) \biggr )\\[6pt]
\nonumber
&=-1+\prod_{k=0}^{n-1}
\biggl (1+\|A\|_{C^0} {\rm Lip}(A) |x-x_0|
\phi(x_0)\phi(x_0+k\alpha) \biggr )\\[6pt]
\nonumber
&\leq -1+\exp \biggl (\sum_{k=0}^{n-1} \|A\|_{C^0} {\rm Lip}(A) |x-x_0|
\phi(x_0)\phi(x_0+k\alpha) \biggr ).
\end{align}
Hence for every $x \in \R/\Z$,
\begin{equation}
\|A_n(x_0)^{-1}(A_n(x)-A_n(x_0))\|\leq
e^{n|x-x_0|\|A\|_{C^0}{\rm Lip}(A)\phi(x_0)S(x_0)}-1,
\end{equation}
which implies
\begin{align}
\|A_n(x)\|&\leq  e^{n|x-x_0|\|A\|_{C^0}{\rm Lip}(A)\phi(x_0)S(x_0)}
\|A_n(x_0)\|\\[6pt]
\nonumber
&\leq  e^{n|x-x_0|\|A\|_{C^0}{\rm Lip}(A)\phi(x_0)S(x_0)} \biggl
(\phi(x_0)\phi(x_0+n\alpha) \biggr )^{1/2}.
\end{align}
\vglue-16pt
\Endproof\vskip4pt 

We now give estimates for the derivatives.

\begin{lemma}
Assume that $A:\R/\Z\to {\rm SL}(2,\R)$ is of class $C^k$ {\rm (}$1 \leq k\break \leq
\infty)$.  Then for every $0\leq r\leq k${\rm ,} and any $x_0,x\in\R/\Z$ such that
$S(x_0)<\infty${\rm ,} 
\begin{equation}
\|(\pa^r A_n)(x)\|\leq C^r n^r \phi(x_0+n\alpha)^{1/2}
\biggl (c_1(x_0)
e^{n c_2(x_0) |x-x_0|} \biggr )^{r+\frac {1} {2}}\|\pa^r
A\|_{C^0}
\end{equation}
where $C$ is an absolute constant and
\begin{align}
c_1(x_0)&=\phi(x_0) S(x_0) \|A\|_{C^0}^2,\\
\nonumber
c_2(x_0)&=2 S(x_0)\phi(x_0)\|A\|_{C^0} \|\pa A\|_{C^0}.
\end{align}
\end{lemma}

%\biggl(C\phi(x_0)S(x_0)
%e^{2n|x-x_0|\phi(x_0)S(x_0) \|A\|_1}\|A\|_0^2\biggr)^{r+1}
%\|A\|_r(\phi(x_0)\phi(x_0+n\alpha))^{1/2}

\Proof 
We compute 
\begin{equation}
\pa^r A_n(x)=\pa^r \left (\prod_{k=n-1}^0 A(\cdot+k\alpha) \right )(x)
\end{equation}
which by Leibniz formula is a sum of $n^r$ terms of the form
\begin{align}
I_{(i^*)}(x)=\left (\prod_{l=n-1}^{i_1+1} A(x+l\alpha) \right ) \cdot
&\pa^{m_1}A(x+i_1\alpha) \cdot
\left (\prod_{l=i_1-1}^{i_2+1} A(x+l\alpha) \right ) \\
\nonumber
&\cdot \pa^{m_2}A(x+i_2\alpha) \cdot
\left (\prod_{l=i_2-1}^{i_3+1} A(x+l\alpha) \right ) \cdots\\
\nonumber
&\cdot \pa^{m_s}A(x+i_s\alpha) \cdot
\left (\prod_{l=i_s-1}^0 A(x+l\alpha) \right )
\end{align}
where $i^*$ runs through ${\cal I}=\{0,\ldots,n-1\}^{\{1,\ldots,r\}}$
and where $s\leq r$ and $\{i_1,\ldots,i_s\}\break =i^*(\{1,...,r\})$
%\mathrm {Im} i^*$,
satisfy $n-1\geq i_1>i_2>\cdots i_s\geq 0$ and
$m_l=\#(i^*)^{-1}(i_l)$. (Notice that $m_1+\ldots+m_s=r$.)
Each term $I_{(i^*)}$ can be written
%\begin{align}
%I_{(i^*)}(x)=\left (\prod_{l=n-1}^0 A(x+l\alpha) \right ) \cdot
%&\Ad \left ( \left (\prod_{l=i_1-1}^0 A(x+l\alpha) \right )^{-1} \right )
%\cdot \biggl(A(x+i_1\alpha)^{-1}\pa^{m_1}A(x+i_1\alpha)\biggr)\cdot\\
%\nonumber
%&\Ad \left ( \left (\prod_{l=i_2-1}^0
%A(x+l\alpha) \right )^{-1} \right )
%\cdot \biggl(A(x+i_2\alpha)^{-1}\pa^{m_2}A(x+i_2\alpha)\biggr)\cdots\\
%\nonumber
%&\Ad \left ( \left (\prod_{l=i_s-1}^0 A(x+l\alpha) \right )^{-1} \right )
%\cdot \biggl(A(x+i_s\alpha)^{-1}\pa^{m_s}A(x+i_s\alpha)\biggr).
%\end{align}
\begin{align}
I_{(i^*)}(x)=A_n(x) \cdot
&\Ad \left (A_{i_1}(x)^{-1} \right )
\cdot \biggl(A(x+i_1\alpha)^{-1}\pa^{m_1}A(x+i_1\alpha)\biggr)\\
\nonumber
&\cdot \Ad \left ( A_{i_2}(x)^{-1} \right )
\cdot \biggl(A(x+i_2\alpha)^{-1}\pa^{m_2}A(x+i_2\alpha)\biggr)\cdots\\
\nonumber
&\cdot \Ad \left ( A_{i_s}(x)^{-1} \right )
\cdot \biggl(A(x+i_s\alpha)^{-1}\pa^{m_s}A(x+i_s\alpha)\biggr).
\end{align}
From the previous lemma,
\begin{eqnarray}
\left \|A_{i_p}(x) \right \|
&\leq&  \bigl(K\phi(x_0)\phi(x_0+i_p\alpha)\bigr)^{1/2},
\\[6pt]
\left \|\Ad \left (A_{i_p}(x) \right ) \right \|
&\leq& \left \|A_{i_p}(x) \right \|^2
\leq K\phi(x_0)\phi(x_0+i_p\alpha)
\end{eqnarray}
where 
\begin{equation}
K=e^{2 n|x-x_0|\phi(x_0)S(x_0)\|A\|_{C^0} \|\partial A\|_{C^0}}.
\end{equation}
Hence we get the following bound
\begin{equation}
\|I_{(i^*)}(x)\|\leq \biggl(K\phi(x_0)\phi(x_0+n\alpha)\biggr)^{1/2}
\prod_{p=1}^s \biggl(
K\phi(x_0)\phi(x_0+i_p \alpha) \|A\|_{C^0} \|\pa^{m_p} A\|_{C^0}\biggr).
\end{equation}
From this and the convexity (Hadamard-Kolmogorov) inequalities \cite {Ko}
\begin{equation}
\|\pa^m A\|_{C^0}\leq C\|A\|_0^{1-(m/r)}\|\pa^r A\|_{C^0}^{\frac {m} {r}},
\quad 0 \leq m \leq r,
\end{equation}
we deduce (using $\sum_{p=1}^s m_p=r$)
\begin{align}
\|I_{(i^*)}(x)\| &\leq \biggl(K \phi(x_0)\phi(x_0+n\alpha)\biggr)^{1/2}\\
&\quad\times
K^s \phi(x_0)^s \|A\|^s_{C^0}\prod_{p=1}^s
\biggl(C \|A\|_{C^0}^{1-\frac {m_p} {r}} \|\pa^r A\|_{C^0}^{\frac
{m_p} {r}} \phi(x_0+i_p \alpha) \biggr)\nonumber\\
\nonumber
&\leq C^s K^{s+\frac {1} {2}} \phi(x_0)^{s+\frac {1} {2}}
\phi(x_0+n\alpha)^{1/2}
\|A\|^{2s-1}_{C^0} \|\pa^r A\|_{C^0} \prod_{p=1}^s \phi(x_0+i_p \alpha)\\
\nonumber
&\leq C^r \bigl(K \|A\|^2_{C^0} \phi(x_0)\bigr)^{r+\frac {1} {2}}
\phi(x_0+n\alpha)^{1/2} \|\pa^r A\|_{C^0} \prod_{p=1}^s
\phi(x_0+i_p\alpha),
\end{align}
so  that
\begin{align}
\|\partial^r A_n(x)\| &\leq \sum_{i^*\in{\cal I}} \|I_{(i^*)}(x)\|\\
\nonumber
&\leq C^r \bigl(K \|A\|^2_{C^0} \phi(x_0) \bigr)^{r+\frac {1} {2}}
\phi(x_0+n\alpha) ^{1/2} \|\pa^r A\|_{C^0}
\\
&\quad\times\nonumber\sum_{i^*\in{\cal I}}\phi(x_0+i_1\alpha)\cdots\phi(x_0+i_s\alpha).
\end{align}
But the last sum in this estimate satisfies the inequality
\begin{multline}
\sum_{i^*\in{\cal I}}\phi(x_0+i_1\alpha)\cdots\phi(x_0+i_s\alpha)
\\ \leq \biggl(\phi(x_0)+\ldots+\phi(x_0+(n-1)\alpha)\biggr)^r \leq
n^r S(x_0)^r
\end{multline}
(recall that $\phi \geq 1$) which implies the \pagebreak result.
\Endproof\vskip4pt 

We can now conclude easily:

\begin{lemma} \label {precompact}

Assume that $A:\R/\Z \to {\rm SL}(2,\R)$ is $C^k$ {\rm (}$1 \leq k \leq \infty${\rm )}.
For almost every $x_* \in \R/\Z${\rm ,} there exists $K>0${\rm ,}
such that for every $d>0$
%\footnote {In what follows, the fact that $K$ is independent
%of $d$ will not be actually used.}
and for every $n>n_0(d)${\rm ,} if $\|\alpha
n\|_{\R/\Z} \leq \frac {d} {n}${\rm ,} then
\begin{equation}
\|\partial^r A_n(x)\| \leq K^{r+1} n^r \|A\|_{C^r},
\quad |x-x_*| \leq \frac {d} {n}.
\end{equation}
\end{lemma}
\vglue8pt

\Proof 
Let $X \subset \R/\Z$ be the set of all $x$ such that $S(x)<\infty$
%$\phi(x)<\infty$
where the $x$ are measurable continuity points of $S$ and
$\phi$.  This means that for every $\epsilon>0$, $x$ is a density
point of
\begin{equation}
Y(x,\epsilon)=S^{-1}(S(x)-\epsilon,S(x)+\epsilon) \cap
\phi^{-1}(\phi(x)-\epsilon,\phi(x)+\epsilon).
\end{equation}
It is a classical fact that $X$ has full Lebesgue Measure.

Fix $x_* \in X$, $d>0$ and $\epsilon>0$.  If $n$ is sufficiently big then
\begin{equation}
\left |Y(x_*,\epsilon) \cap \left [x_*-\frac {2d} {n},x_*+
\frac {2d} {n} \right ] \right | \geq \frac {(4-\epsilon) d} {n}.
\end{equation}
If $\|\alpha n\|_{\R/\Z}<\frac {d} {n}$, this implies
\begin{equation}
\left |(Y(x_*,\epsilon)-\alpha n) \cap Y(x_*,\epsilon)\cap
\left [x_*-\frac {d} {n},x_*+
\frac {d} {n} \right ] \right | \geq \frac {(2-2\epsilon) d} {n}.
\end{equation}
In particular, each point $x \in \left [x_*-\frac {d} {n},x_*+
\frac {d} {n} \right ]$ is at distance at most
$\frac {2 \epsilon d} {n}$ from  a point $x_0$ such that $x_0 \in
Y(x_*,\epsilon)$ and $x_0+\alpha n \in Y(x_*,\epsilon)$.  In particular, for
every $\delta>0$, if $\epsilon>0$ is sufficiently small then
$c_1(x_0) \leq c_1(x_*)+\delta$, $c_2(x_0) \leq c_2(x_*)+\delta$ where
$c_1$ and $c_2$ are as in the previous lemma.  The previous lemma
implies that
\begin{align}
\|(\partial^r A_n)(x)\| &\leq C^r n^r \phi(x_0+n\alpha)^{1/2}
\left (c_1(x_0) e^{c_2(x_0) n |x-x_0|}
\right )^{r+\frac {1} {2}} \|\pa^r A\|_{C^0}\\
\nonumber
&\leq C^r n^r (\phi(x_*)+\epsilon)^{1/2}
\left ((c_1(x_*)+\delta) e^{2 \epsilon d (c_2(x_*)+\delta)}
\right )^{r+\frac {1} {2}} \|\pa^r A\|_{C^0}.
\end{align}

It immediately follows that for every $\epsilon>0$, for every $n$
sufficiently big such that $\|\alpha n\|_{\R/\Z}<\frac {d} {n}$, we have
\begin{equation}
\|\partial^r A_n(x)\| \leq n^r \biggl(C c_1(x_*)+\epsilon
\biggr)^{r+1} \|A\|_{C^r}, \quad |x-x_*| \leq \frac {d} {n}.
\end{equation}
\vglue-24pt
\Endproof\vskip12pt 

\begin{lemma} \label {limit}
Assume that $A:\R/\Z \to {\rm SL}(2,\R)$ is Lipschitz.
For almost every $x_* \in \R/\Z${\rm ,} for every $d>0${\rm ,} for every $\epsilon>0${\rm ,}
if $n>n_0(d,\epsilon)$ and
$\|\alpha n\|_{\R/\Z}\break \leq \frac {d} {n}${\rm ,} then the
matrix $B(x_*) A_n(x) B(x_*)^{-1}$ is $\epsilon$ close to $\SO(2,\R)$
provided that $|x-x_*| \leq \frac {d} {n}$.
\end{lemma}

\Proof 
Let $x_*$ be a measurable continuity point of $S$ and $B$.
By the same argument of the previous lemma, for $n$
big enough, if $\|\alpha n\|_{\R/\Z}<\frac {d} {n}$, then every $x$ such
that $|x-x_*|<\frac {d} {n}$ is at distance at most $\frac {\epsilon} {n}$
from some $x_0$ such that $|S(x_0)-S(x_*)|<\epsilon$,
$\|B(x_0)-B(x_*)\|<\epsilon$ and $\|B(x_0+n\alpha)-B(x_*)\|<\epsilon$.  By
(\ref {lip}), we have
\begin{equation}
\|A_n(x_0)^{-1}(A_n(x)-A_n(x_0))\|\leq
e^{n|x-x_0|\|A\|_{C^0}{\rm Lip}(A)\phi(x_0)S(x_0)}-1 \leq K \epsilon
\end{equation}
and so it is enough to show that $B(x_*) A_n(x_0) B(x_*)^{-1}$ is close to
$\SO(2,\R)$.  But this is clear since $B(x_0+\alpha n) A_n(x_0) B(x_0)^{-1}
\in \SO(2,\R)$ and $B(x_0)$, $B(x_0+n\alpha)$ are close to $B(x_*)$.
\hfq

\section{Renormalization}

Let $\Omega^r=\R \times C^r(\R,{\rm SL}(2,\R))$.  We will view $\Omega^r$ as a
subgroup of\break $\Diff^r(\R \times \R^2)$:
\begin{equation}
(\alpha,A) \cdot (x,w)=(x+\alpha,A(x) \cdot w).
\end{equation}
A $C^r$ fibered $\Z^2$-action is a homomorphism $\Phi:\Z^2 \to
\Omega^r$ (that is,
$\Phi(n,m) \circ \Phi(n',m')=\Phi(n+n',m+m')$).  We let $\Lambda^r$ denote
the space of $C^r$ fibered\break $\Z^2$-actions.  We endow $\Lambda^r$ with the
pointwise topology.  This topology is induced from the embedding
$\Lambda^r \to \Omega^r \times \Omega^r$, $\Phi \mapsto
(\Phi(1,0),\Phi(0,1))$.\footnote{Here and in what follows, spaces of $C^r$
functions (such as $C^r(\R,{\rm SL}(2,\R))$) are always
endowed with the {\it weak} topology of uniform $C^r$-convergence on
compacts.  In the $C^\omega$ case (which is the most important for us), this
means that a sequence $A^{(n)}$ converges to $A$ if and only if for every
compact $K$ there exists a {\it complex} neighborhood $V \supset K$ such
that (the holomorphic extensions of) $A^{(n)}$ (are defined and) converge
to $A$ uniformly on $V$.  We recall that the weak topology is metrizable for
$r \neq \omega$, but not even separable for $r=\omega$.}

Let $$\Pi_1:\R \times C^r(\R,{\rm SL}(2,\R)) \to \R ,\quad 
 \Pi_2:\R \times C^r(\R,{\rm SL}(2,\R)) \to C^r(\R,{\rm SL}(2,\R))$$ be the coordinate
projections.  Let also $\g^\Phi_{n,m}=\Pi_1 \circ \Phi(n,m) \in \R$ and
$A^\Phi_{n,m}=\Pi_2 \circ \Phi(n,m) \in C^r(\R,{\rm SL}(2,\R))$.

The action $\Phi$ will be called nondegenerate
if $\Pi_1 \circ \Phi:\Z^2 \to \R$ is injective.
Let $\Gamma^r$ be the set of nondegenerate actions.

We let $\Lambda^r_0$
be the set of $\Phi \in \Lambda^r$ such that
$\g^\Phi_{1,0}=1$.
For $\Phi \in \Lambda^r_0$, we let $\alpha^\Phi=\g^\Phi_{0,1}$.
We let $\Gamma^r_0=\Gamma^r \cap \Lambda^r_0=\{\Phi \in \Lambda^r_0,\,
\alpha^\Phi \in (0,1) \setminus \Q\}$.

\Subsec{Some operations}
Let $\lambda \neq 0$.
Define $M_\lambda:\Lambda^r \to \Lambda^r$ by
\begin{equation}
M_\lambda(\Phi)(n,m)=(\lambda^{-1} \g^\Phi_{n,m},x \mapsto
A^\Phi_{n,m}(\lambda x)).
\end{equation}

Let $x_* \in \R$.  Define $T_{x_*}:\Lambda^r \to \Lambda^r$ by
\begin{equation}
T_{x_*}(\Phi)(n,m)=(\g^\Phi_{n,m},x \mapsto A^\Phi_{n,m}(x+x_*)).
\end{equation}

Let $U \in \GL(2,\Z)$.  Define $N_U:\Lambda^r \to \Lambda^r$ by
\begin{equation}
N_U(\Phi)(n,m)=\Phi(n',m'), \quad
\begin{pmatrix} n'\\m' \end{pmatrix}=
U^{-1} \cdot
\begin{pmatrix} n\\m \end{pmatrix}.
\end{equation}

The operations $M$, $T$, and $N$ will be called {\it rescaling},
{\it translation}, and {\it base change}.

Notice that $M_\lambda M_{\lambda'}=M_{\lambda \lambda'}$,
$T_{x_*} T_{x'_*}=T_{x_*+x'_*}$,
and $N_U N_{U'}=N_{U U'}$ (that is, $M$, $T$, and $N$ are left actions of
$\R^*$, $\R$ and $\GL(2,\Z)$ on $\Lambda^r$).
Moreover, base changes commute with translations and rescalings.

Notice that $C^r(\R,{\rm SL}(2,\R))$ acts on $\Omega^r$ by
$$\Ad_B(\alpha,A(\cdot))=(\alpha,B(\cdot+\alpha) A(\cdot)
B(\cdot)^{-1}).$$    This action extends to an action (still denoted $\Ad_B$)
on $\Lambda^r$.  We will say that $\Phi$ and $\Ad_B(\Phi)$ are
$C^r$-conjugate via $B$.

\Subsec{Continued fraction expansion}
Let $0<\alpha<1$ be irrational.  We will discuss some elementary facts and
fix notation regarding the continued fraction expansion
\begin{equation} \label {cont}
\alpha=\cfrac {1} {a_1+\cfrac {1} {a_2+\cdots}}
\end{equation}
and we refer the reader to \cite {HW} for details.
Define $\alpha_n=G^n(\alpha)$ where
$G$ is the Gauss map $G(x)=\{x^{-1}\}$ ($\{ \cdot \}$ denotes
the fractionary part).  The coefficients $a_n$ in (\ref {cont}) are given by
$a_n=[\alpha^{-1}_{n-1}]$, where $[ \cdot ]$ denotes the integer part.  We also
set $a_0=0$ for convenience.  Then
\begin{equation}
\alpha_n=\cfrac {1} {a_{n+1}+\cfrac {1} {a_{n+2}+\cdots}}.
\end{equation}

Let $\beta_n=\prod_{j=0}^n \alpha_j$.  Define
\begin{equation}
Q_0=\begin{pmatrix} q_0 & p_0 \\ q_{-1} & p_{-1} \end{pmatrix}=
\begin{pmatrix} 1 & 0 \\ 0 & 1 \end{pmatrix},
\end{equation}
\begin{equation}
Q_n=\begin{pmatrix} q_n & p_n \\ q_{n-1} & p_{n-1} \end{pmatrix}=
\begin{pmatrix} a_n & 1 \\ 1 & 0 \end{pmatrix}
\begin{pmatrix} q_{n-1} & p_{n-1} \\ q_{n-2} & p_{n-2} \end{pmatrix};
\end{equation}
that is,
\begin{equation}
Q_n=U(\alpha_{n-1}) \cdots U(\alpha_0),
\end{equation}
where
\begin{equation} \label {C cdot}
U(x)=\begin{pmatrix} [x^{-1}] & 1\\1 & 0\end{pmatrix}.
\end{equation}
Then we have
\begin{equation}
\beta_n=(-1)^n (q_n \alpha-p_n)=\frac {1} {q_{n+1}+\alpha_{n+1} q_n},
\end{equation}
\begin{equation}
\frac {1} {q_{n+1}+q_n}<\beta_n<\frac {1} {q_{n+1}}.
\end{equation}

\Subsec{Renormalization}
We define the renormalization operator around $0$,
$R \equiv R_0:\Gamma^r_0 \to \Gamma^r_0$, by
$R(\Phi)=M_{\alpha}(N_{U(\alpha)}(\Phi))$ where
$\alpha=\alpha^\Phi$ and $U(\cdot)$ is given by (\ref {C cdot}).

The renormalization operator around $x_* \in \R$, $R_{x_*}:\Gamma^r_0 \to
\Gamma^r_0$ is defined by $R_{x_*}=T_{x_*}^{-1} \circ R \circ T_{x_*}$.

Notice that if $\Phi \in \Gamma^r_0$ and $\alpha^\Phi=\alpha$ then
$\alpha^{R(\Phi)}=G(\alpha)$ and so
\begin{equation}
R^n(\Phi)=M_{\alpha_{n-1}} \circ N_{U(\alpha_{n-1})} \circ \cdots \circ
M_{\alpha_0} \circ N_{U(\alpha_0)}(\Phi)=
M_{\beta_{n-1}}(N_{Q_n}(\Phi)).
\end{equation}

\Subsec{Normalized actions{\rm ,} relation to cocycles}
An action $\Phi \in \Lambda^r_0$ will be called {\it normalized} if
$\Phi(1,0)=(1,\id)$.  If $\Phi$ is normalized then $\Phi(0,1)=(\alpha,A)$
can be viewed as a $C^r$-cocycle, since $A$ is automatically defined modulo
$\Z$.\footnote {Since the commutativity relation
$(1,\id) \circ (\alpha,A)=(\alpha,A) \circ (1,\id)$ is equivalent to
$A(x)=A(x+1)$.}
Inversely, given a $C^r$-cocycle $(\alpha,A)$, $\alpha \in [0,1]$,
we associate a normalized action $\Phi_{\alpha,A}$ by setting
\begin{equation}
\Phi_{\alpha,A}(1,0)=(1,\id), \quad \Phi_{\alpha,A}(0,1)=(\alpha,A).
\end{equation}

\begin{lemma} \label {normalized}
Any $\Phi \in \Lambda^r_0$ is $C^r$-conjugate to a normalized action.
Moreover{\rm ,} if $\Phi_n(1,0) \in \Lambda^r_0$
converges to $(1,\id)$ in $\Lambda^r_0$
then one can choose a sequence of conjugacies
converging to $\id$ in the $C^r$ topology\/\footnote {The reason we refer to
sequences instead of speaking of closeness is because the
$C^\omega$ topology is not separable.}.
\end{lemma}

\Proof 
We first assume that $r \neq \omega$.  Let $\Phi(1,0)=(1,A)$.
Let $B \in C^r([0,3/2],{\rm SL}(2,\R))$
be such that $B(x)=\id$, $x \in [0,1/2]$, $B(x)=A(x-1)$, $x \in [1,3/2]$. 
Let us extend $B$ to $\R$ forcing $\Ad_B (1,\id)=(1,A)$
($B$ is still smooth after the modification).  If $A$ is $C^r$ close to
$\id$, we can select $B:[0,3/2] \to {\rm SL}(2,\R)$ to be $C^r$ close to
$\id$, and in this case $B:\R \to {\rm SL}(2,\R)$ is also $C^r$ close to $\id$.

Let us now assume that $r=\omega$.
Let us first deal with the case where (the holomorphic extension of)
$A$ is close to the identity in a definite neighborhood of $\R$.
Extend $A$ to a real-symmetric \pagebreak
$C^\infty$ function $A:\C \to {\rm SL}(2,\C)$ which is $C^\infty$ close to the
identity and which is holomorphic on a definite
neighborhood $V$ of $\R$.  We will assume that $V$ satisfies (after
shrinking)
\begin{equation}
z \in V \implies z+1 \in V, \quad \Re z \leq 0,
\end{equation}
\begin{equation}
z \in V \implies z-1 \in V,\quad \Re z \geq 1,
\end{equation}
\begin{equation}
[0,1] \times [-\epsilon,\epsilon] \subset V.
\end{equation}

Let $B \in C^\infty(\C,{\rm SL}(2,\C))$ be $C^\infty$ close to the identity,
real-symmetric, and satisfying
$A(z)=B(z+1)B(z)^{-1}$,
$z \in \C$ ($B$ is obtained as in the previous case).
Notice that
$\op B(z+1)=\op A(z)B(z)+A(z)\op B(z)$, so for $z \in V$ we have
$B(z+1)^{-1} \op B(z+1)=B(z+1)^{-1} A(z) \op B(z)=B(z)^{-1} \op B(z)$. 
Moreover,
\begin{equation} \label {Bop}
\|B(z)^{-1} \op B(z)\|<\delta, \quad z \in [0,1] \times
[-\epsilon,\epsilon]
\end{equation}
for some small $\delta$.

Given $C:\R/\Z \times [-1,1] \to {\rm SL}(2,\C)$, we let $D=BC^{-1}$ and
we obviously have $A(z)=D(z+1)D(z)^{-1}$.  We want to choose $C$
so that
\begin{equation} \label {C-1}
\op (C(z)^{-1}) C(z)=-B(z)^{-1} \op B(z), \quad
z \in [0,1] \times [-\epsilon,\epsilon],
\end{equation}
for this will assure us that
\begin{equation}
B(z)^{-1} \op D(z) C(z)=B(z)^{-1} \op B(z)+\op (C(z)^{-1}) C(z)
\end{equation}
vanishes for $z \in [0,1] \times [-\epsilon,\epsilon]$ and also in
$V \cap (\R \times [-\epsilon,\epsilon])$ (this guarantees that $D$ is
holomorphic in a definite neighborhood of $\R$),
and we also want to impose that $C$ (and hence $D$)
is $C^0$ close to the identity.
Here the smoothness requirement on $C$ is for it to be of class $W^{1,1}$;
that is, it should be continuous and have distributional derivatives in
$L^1$.

Equation (\ref {C-1}) is equivalent to
\begin{equation}
C(z)^{-1} \op C(z)=B(z)^{-1} \op B(z).
\end{equation}
To conclude, we use the following proposition:

\begin{prop}
There exists $\kappa>0$ with the following property.  Let $\eta \in
L^\infty(\R/\Z \times [-1,1],\Sl(2,\R))$ and assume that
$\|\eta\|_{L^\infty}<\kappa$.  Then there exists $C:\R/\Z \times
[-1,1] \to {\rm SL}(2,\R)$ of class $W^{1,1}$
such that $C(z)^{-1} \op C(z)=\eta$ and
$\|C-\id\|_{C^0} \leq \kappa^{-1} \|\eta\|_{L^\infty}$ close to the identity
for $z \in \R/\Z \times [-1,1]$.  Moreover{\rm ,} $C$ is real-symmetric
provided $\eta$ is real-symmetric.
\end{prop}

\Proof 
Let $W^{1,1}(\R/\Z \times [-1,1],\Sl(2,\R))$ be the space of continuous
maps $a:\R/\Z \times [-1,1] \to \Sl(2,\R)$
with integrable distributional derivatives, endowed with the natural norm.
We can obtain a bounded linear map
$P:L^\infty(\R/\Z \times [-1,1],\Sl(2,\C)) \to
W^{1,1}(\R/\Z \times [-1,1],\Sl(2,\C))$ 
which is real-symmetric and solves $\op \circ P=\id$.  Indeed $P$ can be
given explicitly in terms
of the Cauchy transform
\begin{equation}
(P \alpha)(z)=\frac {-1} {\pi} \int_{\R \times [-1,1]}
\frac {\alpha(\zeta)} {z-\zeta} d\zeta \wedge
d\overline \zeta=\lim_{t \to \infty}
\frac {-1} {\pi} \int_{[-t,t] \times [-1,1]}
\frac {\alpha(\zeta)} {z-\zeta} d\zeta \wedge
d\overline \zeta.
\end{equation}
Define an analytic map
$T:L^\infty(\R/\Z \times [-1,1]) \to L^\infty(\R/\Z \times [-1,1])$ by
$T(\cdot)=e^{-P(\cdot)} \op e^{P(\cdot)}$.  Then $T(0)=0$, $DT(0)=\id$.
It follows that $T$ is a diffeomorphism in a neighborhood of
$\eta=0$, so we may solve $e^{-P \alpha} \op e^{P \alpha}=\eta$ with
$\|\alpha\|_\infty \leq K \|\eta\|_{L^\infty}$ provided
$\eta$ is close to $0$.
It follows that $C=e^{P \alpha}$ satisfies the conclusion of the
proposition.
\Endproof\vskip4pt 

We may now obtain $C$ with the required properties
by taking $\eta=B^{-1} \op B$ in $[0,1] \times
[-\epsilon,\epsilon]$ and $\eta=0$ otherwise and applying the previous
proposition.  This concludes the second part of the lemma in the case
$r=\omega$.

This argument
also works if we only assume that $A$ is close to the identity
in the $C^\infty$ topology (indeed
the $C^1$ topology is enough, as this is all that we need to get
(\ref {Bop})), and gives the first part of the lemma also in this case
(but we obviously do not get that the holomorphic extension of the
normalizing matrix is close to the identity).
In order to treat the global case, we first consider $B \in
C^\infty(\R,{\rm SL}(2,\R))$ with $A(x)=B(x+1)B(x)^{-1}$, and then
approximate $B$ (in the $C^\infty$ topology) by
$B' \in C^\omega(\R,{\rm SL}(2,\R))$.  Then $B'(x+1)^{-1}A(x)B'(x)$
is $C^\infty$ close to the identity and we can apply the previous case.
\hfq

\Subsec{Degree and fibered rotation number} \label {degree}
The degree and the fibered rotation number of an action will be considered
in detail in Appendix A.  Here we present only a summarized
(and more intuitive) discussion.

The {\it degree} $\deg \Phi$ of a nondegenerate action $\Phi$ can be
defined as follows.  The degree of a normalized
action $\Phi_{\alpha,A}$ is the (topological) degree of the map
$A:\R/\Z \to {\rm SL}(2,\R)$\footnote {Recall that the fundamental group of
${\rm SL}(2,\R)$ is generated by $\theta \mapsto R_\theta$, and hence is
canonically isomorphic to $\Z$.}.  It is easy to see that the
degree of a normalized action is invariant under conjugacies.
This allows us to define the degree of a nondegenerate
action $\Phi$ as the degree of any normalized
action $\Phi_{\alpha,A}$ obtained from $\Phi$ by rescaling and conjugacy.
It is readily seen that the degree is invariant
under rescalings, conjugacies, and translations.  In the Appendix A
we will
see that base changes preserve the degree up to sign:
$\deg N_U(\Phi)=\det U \deg \Phi$.  In particular, the renormalization of an
action of degree $0$ still has  degree $0$.

The {\it fibered rotation number} $\rot(\Phi)$ of an action
$\Phi$ is only defined in the case $\deg \Phi=0$.  For a nondegenerate
action, it can be defined as follows.  If $\Phi$ has degree $0$, and is
conjugated to a normalized action $\Phi_{\alpha,A}$,
then $(\alpha,A)$ is homotopic to the identity, and it is natural to define
$\rot(\Phi)$ as the fibered rotation number of the cocycle $(\alpha,A)$.
In general, a
nondegenerated action $\Phi$ may be rescaled to an action $M_\lambda(\Phi)$
which is conjugated to a normalized action: we then define
$\rot(\Phi)=\lambda \rot(M_\lambda(\Phi))$.
It turns out (see Appendix A) that
$\rot(\Phi)$ is only well defined up to addition of
an element of the {\it module of frequency} of $\Phi$, that is, the
$\Z$-module $\Gamma_\Phi=\{\gamma^\Phi_{n,m},\, (n,m) \in \Z^2\}$, and so
$\rot(\Phi)$ should be regarded as an element of
$\R/\Gamma_\Phi$.\footnote{This is related to the fact that
the fibered rotation number is not a conjugacy invariant for cocycles.}
It is readily seen that
$\rot(\Phi)$ is invariant under translations.  We will see in  Appendix
A
that base changes preserve the fibered rotation number up to sign:
$\rot(N_U(\Phi))=\det U \rot(\Phi)$.

We shall say that an element in $\R/\Gamma_\Phi$ is
{\it Diophantine} if for some representative $\beta$, some basis
$\{e_1,e_2\} \subset \Z^2$ and some
$\kappa>0$, $\tau>0$, one has
\begin{equation}
|2 \beta-k\g_{e_1}-l\g_{e_2}| \geq \kappa (1+|k|+|l|)^{-\tau}, \quad (k,l)
\in \Z^2.
\end{equation}
This definition is clearly independent of the choice of the representative
and of the chosen basis ($\kappa$ then has to be changed).
Finally, we say that the action $\Phi$ is (fiberwise) Diophantine if
$\rot(\Phi)$ is Diophantine.  This notion is stable
under conjugation, translation, rescaling, and base change, so it is also
stable under renormalization.  This definition is such that
a nondegenerate normalized action $\Phi_{\alpha,A}$ is Diophantine if and
only if $\rho(\alpha,A)$ is Diophantine with respect to $\alpha$.

\Subsec{Reducibility}
An action $\Phi$ is called {\it constant} if for every $(n,m) \in \Z^2$, 
$x \mapsto A^\Phi_{n,m}(x)$ is constant.
We will say that an action $\Phi \in \Lambda^r_0$
is {\it $C^r$-reducible} if it is $C^r$-conjugate to a constant action.
It immediately follows that reducibility is invariant
under conjugation, translation, rescaling and base change.  Thus
reducibility is also invariant under renormalization: an action
$\Phi \in \Gamma^r_0$ is\break $C^r$-reducible
if and only if its renormalization $R(\Phi)$ is $C^r$-reducible.
Moreover, reducibility of a {\it nondegenerate
normalized action} $\Phi_{\alpha,A}$
can be interpreted in familiar terms:
%is equivalent to reducibility modulo $\Z$
%of the associated cocycle $(\alpha,A)$.

\begin{lemma}

Let $(\alpha,A) \in (\R \setminus \Q) \times C^r(\R/\Z,{\rm SL}(2,\R))$.  Then
$\Phi_{\alpha,A}$ is $C^r$-reducible
if and only if $(\alpha,A)$ is $C^r$-reducible.

\end{lemma}

\Proof 
Assume that $\Phi_{\alpha,A}$ is reducible.  Then there exists $B\! \in\!
C^r(\R,{\rm SL}(2,\R)\!)$ such that $B(x+1)B(x)^{-1}=U$,
$B(x+\alpha)A(x)B(x)^{-1}=V$, where $U,V \in {\rm SL}(2,\R)$ commute.  Write
$U=\varepsilon e^u$, where $u \in \Sl(2,\R)$ commutes with $V$, and
$\varepsilon \in \{1,-1\}$.
Let $B'(x)=e^{-xu}B(x)$.  Then $B'(x+1)B'(x)^{-1}=\varepsilon \id$, and so
$B'(x+2)=B'(x)$.  Moreover, $B'(x+\alpha)A(x)B'(x)^{-1}=e^{-\alpha u} V$ is
a constant.  Thus $(\alpha,A)$ is reducible.

Assume that $(\alpha,A)$ is reducible.  Thus there exists $B \in
C^r(\R/2\Z,{\rm SL}(2,\R))$ such that $B(x+\alpha)A(x)B(x)^{-1}=C$ for some $C
\in {\rm SL}(2,\R)$.  Let $D(x)=B(x+1)B(x)^{-1}$, so that $D(x+2)=D(x)$.  Then
$CD(x)C^{-1}=D(x+\alpha)$.

Assume that $C$ is not conjugate to a rotation of angle $\theta=k\alpha/2$
for any $k \in \Z \setminus \{0\}$.  Write in the Fourier series
\begin{equation}
D(x)=\sum_{k \in \Z} \hat D(k) e^{\pi i k x}, \quad \hat D(k) \in \mathrm{M}
(2,\C).
\end{equation}
Then
\begin{equation}
\hat D(k) e^{\pi i k \alpha}=C \hat D(k) C^{-1}.
\end{equation}
If $\hat D(k) \neq 0$ for some $k \neq 0$ then $e^{\pi i k \alpha}$ is an
eigenvalue of $\Ad(C):\mathrm{M}(2,\C) \to \mathrm{M}(2,\C)$.  This implies
that $C$ is conjugate to $R_\theta$ where $\theta=\pm \frac {k \alpha} {2}$,
contradicting our assumption.  Thus $D(x)=\hat D(0)$ is a constant, and it
follows that $\Ad_B(\Phi_{\alpha,A})$ is a constant action.

Assume that $C$ is conjugate to a rotation of angle $\theta=k\alpha/2$
for some $k \in \Z \setminus \{0\}$: $C=U R_\theta U^{-1}$, $U \in
{\rm SL}(2,\R)$.  Let $B'(x)=U R_{-(\theta/\alpha)x} U^{-1} B(x)$.  Then
$B'(x+2)=B'(x)$ and $$B'(x+\alpha) A(x) B'(x)^{-1}=U
R_{-(\theta/\alpha)(x+\alpha)} U^{-1} C U R_{(\theta/\alpha)x} U^{-1}=\id.$$ 
Thus, up to changing $B$ to $B'$ we may assume that $C=\id$, and we can
apply the previous case.
\Endproof\vskip4pt 

We will need the following version of a well-known reducibility result:

\begin{lemma} \label {so2r}

Let $\Phi \in \Gamma^r_0$, $r=\omega,\infty$ be $C^r$-conjugate to an
$\SO(2,\R)$ action of degree $0$.
If $\alpha^\Phi \in {\rm DC}$ then $\Phi$ is $C^r$-conjugate to a normalized
constant action.  In particular{\rm ,} $\Phi$ is $C^r$-reducible.

\end{lemma}

\Proof 
We may assume that $\Phi$ is normalized, since we can always conjugate
$\Phi(1,0)$ to $(1,\id)$ via $C^r(\R,\SO(2,\R))$: this
can be done in the same way as in Lemma~\ref {normalized} (it is indeed
easier to proceed for the $\SO(2,\R)$ case).

Let $(\alpha,A)=\Phi(0,1)$, and let $\phi:\R \to \R$ satisfy
$A(x)=R_{\phi(x)}$.
Since $\Phi$ is normalized, $A$ is defined modulo $\Z$,
and since $\Phi$ is of degree
$0$, this implies that $\phi$ is defined modulo $\Z$ as well.

Consider the Fourier series
\begin{equation}
\phi(\theta)=\sum_{k \in \Z} \hat \phi(k) e^{2k\pi i \theta},
\end{equation}
and let
\begin{equation}
\psi(\theta)=\sum_{k \in \Z \setminus \{0\}}
\hat \psi(k) e^{2k\pi i \theta},
\end{equation}
where
\begin{equation}
\hat \psi(k)=\frac {\hat \phi(k)} {1-e^{2 k \pi i \alpha}}, \quad k \neq 0
\end{equation}
so that
\begin{equation}
\phi(x)-\hat \phi(0)=\psi(x)-\psi(x+\alpha).
\end{equation}
The fact that $\alpha \in {\rm DC}$ implies that $|1-e^{2 k \pi i \alpha}|>\kappa
k^{-\tau}$ for some $\kappa>0$, $\tau>0$.  In particular $\psi \in
C^r(\R/\Z,\R)$.

Let $B(x)=R_{\psi(x)}$.
Then $B \in C^r(\R/\Z,\SO(2,\R))$, and we have $B(x+1)\break \cdot B(x)^{-1}=\id$,
$B(x+\alpha)A(x)B(x)^{-1}=R_{\hat \phi(0)}$.  This implies that
$\Ad_B \Phi$ is a normalized constant action.
\Endproof\vskip4pt 

The following is a restatement of Theorem~\ref {eliassoncocycle}
%a result of Eliasson on reducibility of
%cocycles close to constant ones \cite {eliasson}
%(see \cite {AK3} for the $C^\infty$ case)
in the language of actions.

\begin{lemma} \label {eliassonaction}
Let $\Psi \in \Gamma^\omega_0$ be $C^\omega$-conjugate to a normalized
constant action{\rm ,} and let
$\kappa>0$, $\tau>0$ be fixed.  Let $\Psi_n$
be a sequence of Diophantine actions converging to $\Psi$ in
$\Gamma^\omega_0$ and satisfying $\alpha_n \equiv
\alpha^{\Psi_n} \in {\rm DC}(\kappa,\tau)$.
Then $\Psi_n$ is $C^\omega$-reducible for $n$ large enough.
\end{lemma}

\Proof 
After performing a conjugation, we may assume that $\Psi(1,0)=(1,\id)$ and
$\Psi(0,1)=(\hat \alpha,\hat A)$
where $\hat A \in {\rm SL}(2,\R)$ is a constant.  By Lemma~\ref {normalized},
there exists a sequence $B^{(n)} \in C^\omega(\R,{\rm SL}(2,\R))$ converging to
$\id$ which conjugates $\Psi_n$ to a normalized cocycle
$\Psi'_n=\Ad_{B^{(n)}} \Psi_n$.  It follows that
$(\alpha_n,A^{(n)}) \equiv \Psi'_n(0,1)$
converges to $(\hat \alpha,\hat A)$ in the $C^\omega$-topology.
Thus,  Theorem~\ref {eliassoncocycle}
applies and $(\alpha_n,A^{(n)})$ is $C^\omega$-reducible
for $n$ large enough.  This implies that $\Psi'_n$ and $\Psi_n$
are $C^\omega$-reducible as well.
\hfq

\section{{\it A priori\/} bounds and limits of renormalization}

The language of renormalization allows us to restate Lemma~\ref {precompact}
as a precompactness result:

\begin{thm}[{\it A priori\/} bounds] \label {precomp}

Let $\Phi \in \Gamma^r_0$, $r \geq 1${\rm ,}
be a normalized action{\rm ,} and assume that
the cocycle $(\alpha,A)=\Phi(0,1)$ is $L^2$-conjugated to a cocycle
of rotations.  Then for almost every $x_* \in \R${\rm ,} there exists $K>0$ such
that for every $d>0$ and for every $n>n_0(d)${\rm ,}
\begin{equation}
\bigl\|\partial^k A^{R^n_{x_*} \Phi}_{1,0}(x)\bigr\|
\leq K^{k+1} \|A\|_{C^k}, \quad 0 \leq k \leq r, \quad |x-x_*|<d.
\end{equation}
\begin{equation}
\bigl\|\partial^k A^{R^n_{x_*} \Phi}_{0,1}(x)\bigr\|
\leq K^{k+1} \|A\|_{C^k}, \quad 0 \leq k \leq r, \quad |x-x_*|<d.
\end{equation}
In particular{\rm ,} if $r=\omega,\infty$ then $\{R^n_{x_*}(\Phi)\}_n$ is
precompact in $\Lambda^r_0$.
\end{thm}

\Proof 
Apply Lemma~\ref {precompact} to both $(\alpha,A)$ and to $(\alpha,A)^{-1}$,
obtaining a full measure set of ``good points'' $x_*$.
%Let $x_*$ be as in Lemma~\ref {precompact}.
Notice that
\begin{eqnarray}
A^{R^n_{x_*} \Phi}_{1,0}(x)&=&A_{(-1)^{n-1} q_{n-1}}(x_*+\beta_{n-1}(x-x_*)),
\\[6pt]
A^{R^n_{x_*} \Phi}_{0,1}(x)&=&A_{(-1)^n q_n}(x_*+\beta_{n-1}(x-x_*)).
\end{eqnarray}
Fix $d$ (we may assume $d>1$).
Since $\beta_{n-1}<\frac {1} {q_n}<\frac {1} {q_{n-1}}$, the estimates of
Lemma~\ref{precompact} imply that for $0 \leq k \leq r$ and for
$|x-x_*|<d$,
\begin{align}
\bigl\|\partial^k A^{R^n_{x_*} \Phi}_{1,0}(x)\bigr\|
&\leq \beta_{n-1}^k \|(\partial^k A_{(-1)^{n-1} q_{n-1}})
(x_*+\beta_{n-1}(x-x_*))\|\\
\nonumber
&\leq (\beta_{n-1} q_{n-1})^k K^{k+1} \|A\|_{C^k}
\leq K^{k+1} \|A\|_{C^k},
\\[6pt]
\bigl\|\partial^k A^{R^n_{x_*} \Phi}_{0,1}(x)\bigr\|
&\leq \beta_{n-1}^k \|(\partial^k A_{(-1)^n
q_n})(x_*+\beta_{n-1}(x-x_*))\|\\
\nonumber
&\leq (\beta_{n-1} q_n)^k K^{k+1} \|A\|_{C^k} \leq K^{k+1} \|A\|_{C^k}
\end{align}
(notice that $\|A\|_{C^k}=\|A^{-1}\|_{C^k}$).
The precompactness statement is then obvious.
%\footnote{Notice that we do not need to use that $K$ does not
%depend on $d$.}.
\Endproof\vskip4pt 

This result allows us to consider limits of renormalization.  Those are
easy to analyze due to the following simple corollary of Lemma~\ref {limit}:

\begin{thm}[Limits] \label {limits1}

Let $\Phi \in \Gamma^{\mathrm {Lip}}_0$ be a normalized action{\rm ,}
and assume that the cocycle $(\alpha,A)=\Phi(0,1)$ is
$L^2$-conjugated to a cocycle of rotations.
Then for almost every $x_* \in \R${\rm ,} any limit of
$R^n_{x_*}(\Phi)$ is conjugate to an action of rotations{\rm ,}
via a constant $B \in {\rm SL}(2,\R)$.
\end{thm}

We can now prove the following rigidity result.

\begin{thm}[Rigidity] \label {rigidity}
Let $\alpha \in {\rm RDC}${\rm ,} and let $A:\R/\Z \to {\rm SL}(2,\R)$ be $C^\omega$
and homotopic to the identity.
If $(\alpha,A)$ is $L^2$-conjugated to a cocycle of
rotations{\rm ,} and the fibered rotation number of $(\alpha,A)$ is Diophantine
with respect to $\alpha${\rm ,} then $(\alpha,A)$ is $C^\omega$-reducible.
\end{thm}

\Proof 
Let $\alpha \in {\rm RDC}(\kappa,\tau)$ and
let $n_k \to \infty$ be such that $\alpha_{n_k} \in {\rm DC}(\kappa,\tau)$. 
%Let $\Phi=\Phi_{\alpha,A}$.

Consider the renormalizations $\Psi_k=R^{n_k}_{x_*}(\Phi_{\alpha,A})$,
where $x_*$ is as in Theorems \ref {precomp} and \ref {limits1}.
Notice that for every $k$, $\alpha^{\Psi_k} \in {\rm DC}(\kappa,\tau)$ and
$\Psi_k$ is a Diophantine action.
%the fibered rotation number of $\Psi_k$ is Diophantine with respect to
%$\alpha^{\Psi_k}$.

Passing to a subsequence, we may assume that $\Psi_k \to \Psi$ in the
$C^\omega$ topology.  Since ${\rm DC}(\kappa,\tau)$ is compact,
$\alpha^\Psi=\lim \alpha_{n_k} \in {\rm DC}(\kappa,\tau)$.
By Theorem~\ref {limits1}, $\Psi$ is $C^\omega$-conjugate to an $\SO(2,\R)$
action, and so by Lemma~\ref {so2r},
$\Psi$ is\break $C^\omega$-conjugate to a normalized constant action.
Thus Lemma~\ref {eliassonaction} applies
and we conclude that $\Psi_k$ is $C^\omega$-reducible for $k$ large enough.
It follows that $\Phi_{\alpha,A}$ is
reducible, so that $(\alpha,A)$ is reducible as well.
\Endproof\vskip4pt 

  {\it Proof of Theorems {\rm A} and {\rm A}$'$.}
We can now  prove  Theorem A easily.
Let $\alpha \in {\rm RDC}$,
$v \in C^\omega(\R/\Z,\R)$, and let $\Delta$ be the set of $E \in \R$ such that
$(\alpha,S_{v,E})$ is $L^2$-conjugated to a cocycle of
rotations and the fibered rotation number of $(\alpha,S_{v,E})$ is
Diophantine with respect to $\alpha$.  By Theorem \ref {l2} and
Corollary~\ref {dioph2}, $\Delta \cup \{E \in \R,\,
L(\alpha,S_{v,E})>0\}$ has full Lebesgue measure in $\R$,
%By \cite {AS}, \cite {JM},
%see Appendix,
%$N(E)$ coincides with
%the fibered rotation number of $(\alpha,A)$, so
and Theorem~\ref {rigidity}
implies that $(\alpha,S_{v,E})$ is $C^\omega$-reducible for all $E \in \Delta$.
This shows that $(\alpha,S_{v,E})$ is $C^\omega$-reducible for almost every $E
\in \R$ such that $L(\alpha,S_{v,E})=0$.
By Remark~1.3, if $E \in \R$ is such that
$L(\alpha,S_{v,E})>0$ then $(\alpha,S_{v,E})$
is either nonuniformly hyperbolic or $C^\omega$-reducible,
and the result follows.

This argument also works for Theorem A$'$, if we use Theorem \ref {l3} and
Corollary \ref {dioph3} instead of Theorem~\ref {l2} and
Corollary \ref {dioph2}.
\hfq

\demo{Acknowledgements} We would like to thank Hakan Eliasson, Svetlana
Jitomirskaya, Barry Simon, and Jean-Christophe Yoccoz for several
discussions and suggestions.  We also thank the referee, whose comments
were useful in improving the presentation of this article.

\appendix

\vglue16pt\centerline{\bf Appendix A. Degree and and fibered rotation number} 
\vglue12pt
\setcounter{section}{1}
\setcounter{equation}{0}

In this section we will recall the intrinsic definition of degree and
fibered rotation number for actions given in \cite {K2}, and check that they
coincide with the definitions given in \S \ref {degree}.
The advantage of the
intrinsic definitions is that they allow us to compute easily the effect of
base changes.

For $\alpha \in \R$ and $A:\R \to {\rm SL}(2,\R)$ continuous, we introduce the
following objects.
If $w$ is a point of the  usual euclidean circle 
$\Sr^1 \subset \R^2 \equiv \C$ we set
\begin{equation}
f^A(x,w)=\frac{A \cdot w}
{\|A \cdot w\|},
\end{equation}
and  define, for $\alpha \in \R$,
\begin{align}
F^{\alpha,A}:\R\times\Sr^1&\to {\R}\times \Sr^1\\
\nonumber
(x,w)&\mapsto (x+\alpha,f^A(x,w)).
\end{align}  
If $\pi:\R \to \Sr^1$ is the projection $\pi(y)=\exp(2\pi i y)$ we can
find a continuous lift $d^A:\R \times \R \to \R$ of
$f^A(x,w) w^{-1}$, that is
\begin{equation}
\pi(y+d^A(x,y))=f^A(x,\pi(y)).
\end{equation}
Observe that such a lift is not uniquely defined, every other lift being of
the form $d^A(x,y)+k$, where $k$ is a
constant integer.  Also, for any $x,y \in \R\times \R$ we have
$d^A(x,y+1)=d^A(x,y)$
and thus $d^A(x,w)$ can be defined for any $x \in \R$,
$w \in \Sr^1$.

\Subsec{Cocycles} \label {A1}
Let us first consider the case of a cocycle $(\alpha,A) \in \R/\Z \times
C^0(\R/\Z,{\rm SL}(2,\R))$.  Viewing $A$ as defined on $\R$, we can define $d^A$
(up to an integer), and we get $d^A(x+1,w)=d^A(x,w)+n$, where $n$ is the
topological degree of $\R/\Z \to {\rm SL}(2,\R)$.  Indeed, up to homotopy, we may
assume that $A(x)=R_{n x}$, and we have $d^A(x,w)=nx$.

If $(\alpha,A)$ is homotopic to the identity, $d^A$ descends
to a map $\R/\Z \times \Sr^1 \to \R$ and $F^{\alpha,A}$ descends
to a map $\R/\Z \times \Sr^1 \to \R/\Z \times \Sr^1$.
The usual definition (see \cite {H},
\cite {JM}) of the fibered
rotation number of $(\alpha,A)$ is
\begin{equation}
\rho(\alpha,A)=\int_{\R/\Z \times \Sr^1} d^A(x,w) d\mu(x,w),
\end{equation}
(defined modulo an integer)
where $\mu$ is any probability measure which is invariant under $F^A:\R/\Z
\times \Sr^1 \to \R/\Z \times \Sr^1$ and which projects to Lebesgue measure
on $\Sr^1$.  One easily checks that if $x \in \R$ and $w,w' \in \Sr^1$ then
$|\sum_{k=0}^{n-1} d^A \circ (F^{\alpha,A})^k(x,w')-
\sum_{k=0}^{n-1} d^A \circ (F^{\alpha,A})^k(x,w)|<1$.  This implies that
\begin{align}
&\left |n \int_{\R/\Z \times \Sr^1} d^A(x,w) d\mu(x,w)-
\int_{\R/\Z \times \Sr^1} \sum_{k=0}^{n-1} d^A \circ (F^{\alpha,A})^k(x,w)
d\mathrm{Leb}(x,w) \right |\\
\nonumber
&=\left |\int_{\R/\Z \times \Sr^1} \sum_{k=0}^{n-1} d^A \circ
(F^{\alpha,A})^k(x,w) d\mu(x,w)\right.\\
&\qquad \left.\nonumber -\int_{\R/\Z \times \Sr^1}
\sum_{k=0}^{n-1} d^A \circ (F^{\alpha,A})^k(x,w) d\mathrm{Leb}(x,w)\right |<1,
\end{align}
for every $n\!>\!0$, so that $\rho(\alpha,A)\!=\!\lim \frac {1} {n}
\int_{\R/\Z \times \Sr^1} \sum_{k=0}^{n-1} d^A \circ (F^{\alpha,A})^k(x,w)
d\mathrm{Leb}(x,w)$ does not depend on $\mu$.

\Subsec{Actions}
 Let $(e_1,e_2)$ be a basis of the $\Z$-module $\Z^2$.
Then it is easy to see that the quantity
\begin{equation} \label {e1 e2}
\deg_{e_1,e_2} \Phi=
(d^{A^\Phi_{e_1}}\circ F^{\Phi(e_2)}+d^{A^\Phi_{e_2}})-(d^{A^\Phi_{e_2}}
\circ F^{\Phi(e_1)}+d^{A^\Phi_{e_1}})
\end{equation}
is independent of the choices made for the lifts
and is a constant integer.  Obviously from (\ref {e1 e2}),
$\deg_{e_2,e_1} \Phi=-\deg_{e_1,e_2} \Phi$.
Notice that $d^{A^\Phi_{e_1+e_2}}=
d^{A^\Phi_{e_1}} \circ F^{\Phi(e_2)}+d^{A^\Phi_{e_2}}$
(up to a constant integer), so that
\begin{align}
&\deg_{e_1,e_1+e_2} \Phi=(d^{A^\Phi_{e_1}}\circ
F^{\Phi(e_1+e_2)}+d^{A^\Phi_{e_1+e_2}})-(d^{A^\Phi_{e_1+e_2}}\circ
F^{\Phi(e_1)}+d^{A^\Phi_{e_1}})\\
\nonumber
&=(d^{A^\Phi_{e_1}}\circ
F^{\Phi(e_1+e_2)}+d^{A^\Phi_{e_1}} \circ F^{\Phi(e_2)}+d^{A^\Phi_{e_2}})\\
\nonumber
&\quad -
(d^{A^\Phi_{e_1}} \circ F^{\Phi(e_2)} \circ F^{\Phi(e_1)}+d^{A^\Phi_{e_2}}
\circ F^{\Phi(e_1)}+d^{A^\Phi_{e_1}})\\
\nonumber
&=(d^{A^\Phi_{e_1}}\circ
F^{\Phi(e_1+e_2)}+d^{A^\Phi_{e_1}} \circ F^{\Phi(e_2)}+d^{A^\Phi_{e_2}})\\
\nonumber
&\quad -
(d^{A^\Phi_{e_1}} \circ F^{\Phi(e_1+e_2)}+d^{A^\Phi_{e_2}}
\circ F^{\Phi(e_1)}+d^{A^\Phi_{e_1}})\\
\nonumber
&=(d^{A^\Phi_{e_1}} \circ F^{\Phi(e_2)}+d^{A^\Phi_{e_2}})-
(d^{A^\Phi_{e_2}} \circ F^{\Phi(e_1)}+d^{A^\Phi_{e_1}})=\deg_{e_1,e_2}(\Phi).
\end{align}
A similar computation gives $\deg_{e_1,-e_2} \Phi=-\deg_{e_1,e_2} \Phi$.
These  elementary base change
rules imply that $\deg_{U \cdot e_1,U \cdot e_2} \Phi=\det U
\deg_{e_1,e_2} \Phi$ for any $U \in \GL(2,\Z)$.

We define  $\deg \Phi$ as
$\deg_{(0,1),(1,0)} \Phi$.  To see that this coincides with the
previous definition (given in \S \ref {degree}),
it is enough to check it in the case of a
normalized action $\Phi=\Phi_{\alpha,A}$.  Recalling that
$d^A(x+1,w)=d^A(x,w)+n$ where $n$ is the topological degree of $A:\R/\Z \to
{\rm SL}(2,\R)$, we get from
$d^{A^\Phi_{(1,0)}}=0$ and $d^{A^\Phi_{(0,1)}}(x,w)=d^A(x,w)$ that
$\deg(\Phi)=d^A(x+1,w)-d^A(x,w)=n$,
according to the previous definition.

Assume now that the action $\Phi$ has {\it degree zero}.  Let us denote by
${\cal M}^\Phi$ the set of measures on $\R\times\Sr^1$ which  project
on the first factor to Lebesgue measure on $\R$ and which are invariant
by $F^{\Phi(n,m)}$ for any
$(n,m)\in\Z^2$.  It is not difficult to see that $\MM^\Phi$ is nonempty.
Take as before $(e_1,e_2)$ to be a basis of
$\Z^2$, and for $\mu \in \MM^\Phi$, define the quantity:
\begin{equation}
\rot_{e_1,e_2,\mu} \Phi=
I(0,\g_{e_2}^\Phi;d^{A^\Phi_{e_1}},\mu)-I(0,\g_{e_1}^\Phi;
d^{A^\Phi_{e_2}},\mu),
\end{equation}
where we have defined for any function
$h:\R\times\Sr^1 \to \R$ and $(a,b)\in \R^2$ the quantity
\begin{equation}
I(a,b;h,\mu)={\rm sgn}(b-a)\int_{[a,b]\times\Sr^1}h(x,v)d\mu(x,v).
\end{equation}
If we make other choices for the lifts of $F^\Phi$, the numbers we obtain
just differ by the addition of an element of the
{\it module of frequency} of $\Phi$.

We notice that $\rot_{e_2,e_1,\mu} \Phi=-\rot_{e_1,e_2,\mu} \Phi$ and
\begin{align}
&\rot_{e_1,e_1+e_2,\mu}
\Phi=I(0,\g^\Phi_{e_1}+\g^\Phi_{e_2};d^{A_{e_1}^\Phi},\mu)-
I(0,\g^\Phi_{e_1};d^{A^\Phi_{e_1+e_2}},\mu)\\
\nonumber
&=I(0,\g^\Phi_{e_2};d^{A_{e_1}^\Phi},\mu)+I(\g^\Phi_{e_2},\g^\Phi_{e_1}+
\g^\Phi_{e_2};
d^{A_{e_1}^\Phi},\mu)-I(0,\g^\Phi_{e_1};d^{A^\Phi_{e_1}} \circ
F^{\Phi(e_2)}+d^{A^\Phi_{e_2}},\mu)\\
\nonumber
&=\rot_{e_1,e_2,\mu} \Phi+
I(\g^\Phi_{e_2},\g^\Phi_{e_1}+\g^\Phi_{e_2};   
d^{A_{e_1}^\Phi}\!,\mu)-I(0,\g^\Phi_{e_1};d^{A^\Phi_{e_1}} \circ
F^{\Phi(e_2)}\!,\mu)=\rot_{e_1,e_2,\mu} \Phi,
\end{align}
since
\begin{eqnarray}
\int_{[0,\g^\Phi_{e_1}] \times \Sr^1} d^{A^\Phi_{e_1}} \circ F^{\Phi(e_2)}
d\mu&=&
\int_{F^{\Phi(e_2)}([0,\g^\Phi_{e_1}] \times \Sr^1)} d^{A^\Phi_{e_1}}
d(F^{\Phi(e_2)})_* \mu\\
\nonumber &=&\int_{[\g^\Phi_{e_2},\g^\Phi_{e_1}+\g^\Phi_{e_2}]
\times \Sr^1} d^{A^\Phi_{e_1}} d\mu.
\end{eqnarray}
A similar computation gives $\rot_{e_1,-e_2,\mu} \Phi=-\rot_{e_1,e_2,\mu} \Phi$.
Those elementary base change rules imply that
$\rot_{U \cdot e_1,U \cdot e_2,\mu} \Phi=
\det U \rot_{e_1,e_2,\mu} \Phi$ for any $U \in \GL(2,\Z)$.

Given $B:\R \to {\rm SL}(2,\R)$ continuous, we notice that $F^{0,B}_* \MM^{\Phi}=
\MM^{\Ad_B \Phi}$, and it follows immediately from the definition that
$$\rot_{e_1,e_2,\mu} \Phi=\rot_{e_1,e_2,F^{0,B}_* \mu} \Ad_B \Phi.$$  The
transformation rule for $M_\lambda$ can be also readily checked: $\rot
M_\lambda(\Phi)=\lambda^{-1} \rot \Phi$.

Let us check that $\rot_{e_1,e_2,\mu} \Phi$ does not depend on
$\mu \in \MM$.  This is obvious if $\Gamma_\Phi=\{0\}$ (in this case
$\rot=0$).  Otherwise, via conjugacies, scalings, and base change, we reduce
to the case of checking that $\rot_{(0,1),(1,0),\mu} \Phi$ does not depend
on $\mu$ when $\Phi$ is a normalized action $\Phi_{\alpha,A}$.  In this
case, measures in $\MM^\Phi$ are invariant under $(x,w) \mapsto (x+1,w)$, and so
they descend to $\R/\Z \times \Sr^1$.  Since $A:\R/\Z \to {\rm SL}(2,\R)$ is
homotopic to the identity, we have $d^A(x+1,w)=d^A(x,w)$, so that  $d^A$ also
descends to $\R/\Z \times \Sr^1$.  We have
\begin{equation}
\rot_{(0,1),(1,0),\mu} \Phi=I(0,1;d^A,\mu)=
\int_{\R/\Z \times \Sr^1} d^A(x,w) d\mu(x,w).
\end{equation}
This is precisely the usual definition of the fibered rotation number
$\rho(\alpha,A)$ (see \S \ref {A1}), which does not depend on
$\mu$.  This also shows that setting $\rot \Phi=\rot_{(0,1),(1,0),\mu} \Phi$
one recovers the previous definition (given in \S \ref {degree})
of the fibered rotation number of a
nondegenerate action.

\end{document}